\documentclass[hidelinks,onefignum,onetabnum]{siamart251216}

\usepackage{amsfonts,mathrsfs}
\usepackage{graphicx}
\usepackage{epstopdf}
\usepackage{algorithmic}
\usepackage{tikz}
\usetikzlibrary{calc,positioning,shapes.geometric,arrows.meta,decorations.markings}

\ifpdf
  \DeclareGraphicsExtensions{.eps,.pdf,.png,.jpg}
\else
  \DeclareGraphicsExtensions{.eps}
\fi

\newsiamremark{remark}{Remark}
\newsiamremark{hypothesis}{Hypothesis}
\crefname{hypothesis}{Hypothesis}{Hypotheses}
\newsiamthm{claim}{Claim}
\newsiamremark{fact}{Fact}
\crefname{fact}{Fact}{Facts}

\headers{Energy Dissipation for Geometric Flows}
{Z. Pan, W. Jiang, and J. Jia}

\title{Energy-Dissipative, High-Order BDF Schemes for Simulating Geometric Flows via Energy Increments}

\author{
	Zhiqing Pan\thanks{
		Department of Computational Mathematics, School of Mathematics,
		Jilin University, Changchun, China
		(\email{panzq24@mails.jlu.edu.cn}).
	}
	\and
	Wei Jiang\thanks{
		School of Mathematics and Statistics, Hubei Key Laboratory of Computational Science,
		Wuhan University, Wuhan 430072, China
		(\email{jiangwei1007@whu.edu.cn}).
	}
	\and
	Jiwei Jia\thanks{Corresponding author. Department of Computational Mathematics, School of Mathematics,
		Jilin University, Changchun, China; AI for Science and Engineering Center,
		Shenzhen Loop Area Institute, Shenzhen 518048, China
		(\email{jiajiwei@jlu.edu.cn}).
	}
}

\usepackage{amsopn}

\ifpdf
\hypersetup{
  pdftitle={Energy-dissipative, high-order BDF schemes for simulating geometric flows via energy increments},
  pdfauthor={Zhiqing Pan, Wei Jiang, Jiwei Jia}
}
\fi

\begin{document}
	
\maketitle
	
\begin{abstract}
We propose energy-dissipative backward differentiation formula (BDF)
methods for geometric flows in the Barrett--Garcke--N\"urnberg (BGN) parametric
finite element framework. Starting from a Lagrange-multiplier formulation,
we replace the multistep approximation of the energy derivative by an
energy-increment equation while retaining the BDF geometric equations.
The Simpson quadrature and two interleaved time grids provide third- and
fourth-order constructions. The positive quadrature weights ensure
stepwise decay of the original discrete geometric energy for every
nondegenerate solution of the nonlinear system. We establish an explicit
conditional consistency bound for the scalar energy equation and a
regularity criterion for the Newton system of curve-shortening flows. 
Then, the proposed schemes extend to simulating
curve diffusion, mean curvature flow, and
surface diffusion. Finally, numerical experiments demonstrate geometric accuracy,
energy decay, good mesh distribution, and high performance of the proposed schemes.
\end{abstract}



\section{Introduction}

Curvature-driven flows provide models for interface motion in materials
science, biomembrane dynamics, and image processing
\cite{cao2003geometric,davi1990motion,elliott2010modeling}. Many of these
models are gradient flows of a geometric energy $W$ and satisfy
\begin{equation}\label{eq:intro-energy-law}
\frac{d}{dt}W(\Gamma(t))=-\mathcal D(\Gamma(t)),\qquad
\mathcal D(\Gamma(t))\geq0,
\end{equation}
where $\Gamma(t)$ is the evolving curve or surface and $\mathcal D$ is the
dissipation rate. Curve shortening and curve diffusion flows, for example,
decrease perimeter with dissipation rates $\int_\Gamma\kappa^2\,ds$ and
$\int_\Gamma|\partial_s\kappa|^2\,ds$, respectively; analogous laws hold
for mean curvature flow and surface diffusion
\cite{barrett2020parametric}. A numerical method should reproduce this
energy decrease while resolving the interface motion accurately.
Combining higher-order time integration with discrete energy dissipation
and a well-distributed moving mesh is the problem addressed here.

Parametric finite element methods (PFEMs) approximate the interface
directly by a polygonal curve or a triangulated surface. Their
performance depends on how the mesh evolves along with the geometry.
The Barrett--Garcke--N\"urnberg (BGN) formulation
\cite{barrett2007parametric,barrett2008parametric1,barrett2008parametric}
couples the normal-velocity equation to a weak curvature relation.
The latter also controls tangential nodal motion, which influences the
discrete mesh although it leaves the smooth evolving shape unchanged.
This mechanism underlies the favorable mesh redistribution of classical
BGN schemes. Alternative approaches prescribe artificial tangential
motion through harmonic-map or minimal-deformation formulations
\cite{duan2024new,gao2026energy,hu2022evolving}.
Structure-preserving PFEMs have also been developed for isotropic and
anisotropic surface diffusion
\cite{Bao2023Symmetrized,bao2025unified,bao2021structure}.

Higher-order BGN methods must approximate both the interface velocity and
the geometry on which the discrete equations are evaluated. Jiang, Su,
and Zhang developed a second-order method with interlaced perimeter decay
\cite{jiang2024second} and backward differentiation formula (BDF) methods
of high orders from the second to fourth \cite{jiang2024stable}. In the latter,
lower-order BGN solves supply the predicted integration geometry.
A related predictor--corrector construction was developed for curve and
surface diffusion \cite{jiang2025predictor}. These developments motivate
retaining the BGN spatial structure while improving the treatment of
its energy dissipation law.

Lagrange multipliers provide a direct way to impose geometric balance
laws as additional equations. Garcke, Jiang, Su, and Zhang
\cite{garcke2025structure} used this approach to enforce energy
dissipation and area or volume conservation in PFEMs for surface
diffusion. Their BDF formulation applies the multistep difference
operator to the energy itself. For BDF2, a nonincreasing starting step
allows an induction argument for subsequent energy decay. At orders
three and four, however, a nonpositive BDF approximation of the energy
derivative does not by itself imply $W^{m+1}\leq W^m$.
The issue is the sign of the
individual energy increments in the multistep recurrence, rather than
the order of approximation of the derivative. Moreover, replacing this
equation by a backward-Euler energy update would generally reduce the
consistency order of the coupled method.

Several other approaches achieve high-order energy stability for phase-field gradient flow models. 
On fixed domains, scalar auxiliary variable (SAV) methods
\cite{huang2022new,shen2018scalar,shen2019new} and energy-diminishing
Runge--Kutta methods \cite{fu2025higher,fu2024energy} provide important
examples. For evolving interfaces, Guo and Li \cite{guo2025structure}
combine BDF integration with minimal-deformation, SAV, and
Lagrange-multiplier formulations for anisotropic surface diffusion.
The dissipation theorem for their SAV variants concerns an auxiliary
energy, which is distinct from monotonicity of the original discrete
geometric energy.

Duan, Li, and Zhang \cite{duan2021highorder} obtained high-order
energy-diminishing methods for curvature-dependent surface energies by
averaged vector-field collocation. More recently, Ma and Rao
\cite{ma2026highorder} constructed algebraically stable implicit
Runge--Kutta extensions of BGN-type methods for curves and genus-zero
surfaces. Zhang, Andrews, and Farrell \cite{zhang2026arbitrary} developed
continuous Petrov--Galerkin time discretizations of arbitrary order
with auxiliary variables that preserve energy
dissipation and, for surface diffusion, enclosed volume of the surface. These methods establish
that high-order geometric energy dissipation is achievable. In this paper, our objective
is its realization within the BDF3 and BDF4 Lagrange-multiplier
formulations, with their geometric evolution and curvature equations
retained.

We replace the multistep energy-derivative equation by an approximation
of the integrated balance in \eqref{eq:intro-energy-law}. Writing
$W^m=W(\Gamma^m)$ and $\tau=t_{m+1}-t_m$, the scalar update takes the form
\begin{equation}\label{eq:intro-increment}
W^{m+1}-W^m=-\tau\mathcal Q_k^{m+1},
\end{equation}
where $\mathcal Q_k^{m+1}$ approximates the time-averaged dissipation rate by
a quadrature rule with positive weights. For $k=3,4$, Simpson's rule
combines endpoint and midpoint dissipation values. An auxiliary
half-integer time grid supplies the midpoint states, and the two grids
are advanced alternately. The energy in \eqref{eq:intro-increment} is, respectively,
the length for curve evolution, area for surface evolution, or anisotropic free energy for anisotropic geometric flows.
Thus, every nondegenerate solution satisfying the scalar equation
inherits the corresponding original energy decrease rate, separately on each time
level. The construction only changes the scalar balance rather than the BDF
approximation of the geometric evolution.

The analysis identifies the conditions needed to interpret this
construction. We establish stepwise energy dissipation and quantify the
scalar consistency defect in terms of the quadrature error and the
errors in the nodal dissipation data. We also derive the Newton system
for curve shortening flow (CSF), give a Schur-complement criterion for
its local regularity, and relate the nonlinear energy residual to the
computed energy balance. These results make explicit the roles of
predicted geometry, midpoint accuracy, initialization, and nonlinear
solution accuracy.

The rest of the paper is organized as follows. Section~\ref{sec:method} introduces the energy-increment BDF schemes and
their scalar analysis. Section~\ref{sec:staggered} presents the
interleaved implementation, and Section~\ref{sec:newton} develops the
Newton solver. Extensions to other geometric flows are given in
Section~\ref{sec:extensions}, followed by numerical experiments in
Section~\ref{sec:numerics} and conclusions in
Section~\ref{sec:conclusions}.

\section{Energy-increment BDF discretization}\label{sec:method}

We use curve-shortening flow to introduce the construction. After fixing
the geometric conventions and spatial discretization, we examine the
energy constraint in derivative-based BDF schemes and replace it by an
integrated balance. The resulting scalar equation admits separate
consistency and dissipation arguments.

\subsection{Curve shortening flow and spatial semidiscretization}

Let $\Gamma(t)\subset\mathbb R^2$ be a smooth, simple closed curve with
initial configuration $\Gamma(0)=\Gamma_0$ and regular parametrization
\[
	\boldsymbol{X}(\rho,t):=(x_1(\rho,t),x_2(\rho,t)):
	\mathbb{I}\times[0,T]\to\mathbb{R}^2,
	\qquad \mathbb{I}=[0,1]\ \text{with endpoints identified}.
\]
We orient the curve counterclockwise and choose the outward unit normal
$\mathbf n$. The convention $\kappa\mathbf n=-\partial_{ss}\mathbf X$
then assigns curvature $\kappa=1/R$ to a circle of radius $R$.
All continuous identities are understood while the curve remains smooth
and regular. The arc-length variable is
$
s(\rho,t):=\int_0^\rho |\partial_q \boldsymbol{X}|\,dq
$,
with derivative and measure
$
\partial_s:=\frac{1}{|\partial_\rho \boldsymbol{X}|}\partial_\rho,
\quad
ds=|\partial_\rho \boldsymbol{X}|\,d\rho.
$

Write $W(t)=L(\Gamma(t))$ for the perimeter and
$V=\partial_t\mathbf X\cdot\mathbf n$ for the normal velocity.
The perimeter variation underlying the Lagrange-multiplier approach
\cite{garcke2025structure} is
\[
	\frac{d}{dt}L(\Gamma(t))
	=
	\int_{\Gamma(t)}\kappa V\,ds.
\]
Introducing a scalar multiplier $\lambda(t)$ in the curvature direction
and imposing the perimeter dissipation law gives
\begin{subequations}\label{lagrange multiplier general form}
	\begin{align}
		\partial_t \mathbf{X}(\cdot, t) \cdot \mathbf{n}
		&= -\kappa + \lambda(t)\kappa,
		\label{eq:continuous-evolution} \\
		\kappa \mathbf{n} &= -\partial_{ss} \mathbf{X},
		\label{eq:continuous-geometric} \\
		\frac{dL}{dt} &= -\int_{\Gamma(t)}\kappa^2\,ds,
		\label{eq:continuous-energy-dissipation}
	\end{align}
\end{subequations}
Combining the perimeter variation with
\eqref{eq:continuous-evolution} and
\eqref{eq:continuous-energy-dissipation} yields
$\lambda(t)\int_{\Gamma(t)}\kappa^2\,ds=0$.
Thus $\lambda(t)=0$ whenever $\kappa\not\equiv0$, and
\eqref{lagrange multiplier general form} recovers classical
curve-shortening flow.

For the spatial discretization, partition $\mathbb I$ uniformly into
\[
	\mathbb{I}=[0,1]=\bigcup_{j=1}^{N} I_j,
	\qquad
	I_j=[\rho_{j-1},\rho_j],\qquad \rho_j=jh,\quad h=\frac1N.
\]
Define the finite element space
\[
	V^h:=\left\{u\in C(\mathbb{I}):
	u|_{I_j}\in \mathbb{P}_1,\ \forall j=1,\dots,N,\
	u(\rho_0)=u(\rho_N)\right\}\subseteq H^1(\mathbb{I}),
\]
where $\mathbb{P}_1$ denotes the set of polynomials of degree at most one.
For $\boldsymbol X^h\in[V^h]^2$, let
$\Gamma^h=\boldsymbol X^h(\mathbb I)$ be the associated polygon, with edges
$\boldsymbol{h}_j=\boldsymbol{X}^h(\rho_j)-\boldsymbol{X}^h(\rho_{j-1})$.
All polygons used to evaluate inner products or arc-length derivatives
are assumed to have positive edge lengths. The discrete perimeter is
\begin{equation}\label{eq:discrete-perimeter}
	W^h(t):=L(\Gamma^h(t))
	:=\int_{\Gamma^h(t)}1\,ds
	=\sum_{j=1}^{N}|\mathbf h_j(t)|.
\end{equation}
On $I_j$, the discrete arc-length derivative is
$\partial_s=(h/|\mathbf h_j|)\partial_\rho$.
We use the standard inner product
$(u,v)_{\Gamma^h}:=\int_{\Gamma^h}u\cdot v\,ds$
and its mass-lumped counterpart
\begin{equation}
	(u,v)_{\Gamma^h}^h
	:=
	\frac12 \sum_{j=1}^{N} |\boldsymbol{h}_j|
	\left[(u\cdot v)(\rho_j^-)+(u\cdot v)(\rho_{j-1}^+)\right],
\end{equation}
for scalar- or vector-valued piecewise continuous functions.
The values at $\rho_j^\pm$ denote one-sided traces, with periodic
extension at $\rho_0$ and $\rho_N$.

Given compatible initial data, seek
$
(\boldsymbol{X}^h,\kappa^h,\lambda^h)
\in [V^h]^2\times V^h\times \mathbb{R}
$
such that, for all $(\varphi^h,\boldsymbol{\omega}^h)
\in V^h\times [V^h]^2$,
\begin{subequations}\label{eq:semi}
	\begin{align}
		&\left(\partial_t \boldsymbol{X}^h,
		\varphi^h \boldsymbol{n}^h\right)_{\Gamma^h}^h
		+(1-\lambda^h)\left(\kappa^h,\varphi^h\right)_{\Gamma^h}^h=0,
		\label{eq:semia}\\
		&\left(\kappa^h \boldsymbol{n}^h,
		\boldsymbol{\omega}^h\right)_{\Gamma^h}^h
		-\left(\partial_s \boldsymbol{X}^h,
		\partial_s \boldsymbol{\omega}^h\right)_{\Gamma^h}=0,
		\label{eq:semib}\\
		&\frac{dW^h}{dt}
		+\left(\kappa^h,\kappa^h\right)^h_{\Gamma^h}=0.
		\label{eq:semic}
	\end{align}
\end{subequations}
Here $\Gamma^h(t)=\boldsymbol X^h(\mathbb I,t)$, and the elementwise
outward unit normal is
\[
	\left.\boldsymbol{n}^h\right|_{I_j}
	=
	-\frac{\boldsymbol{h}_j^\perp}{|\boldsymbol{h}_j|},
	\qquad j=1,\dots,N,
\]
with $(v_1,v_2)^\perp=(-v_2,v_1)$ denoting counterclockwise rotation by
$\pi/2$. This sign is consistent with the counterclockwise orientation.

For a differentiable polygonal evolution,
$dW^h/dt=(\partial_s\mathbf X^h,\partial_s\partial_t\mathbf X^h)_{\Gamma^h}$.
Testing \eqref{eq:semia} with $\kappa^h$ and \eqref{eq:semib} with
$\partial_t\mathbf X^h$ therefore gives
$dW^h/dt=-(1-\lambda^h)(\kappa^h,\kappa^h)^h_{\Gamma^h}$.
Equation~\eqref{eq:semic} then implies $\lambda^h=0$ when
$\kappa^h\not\equiv0$. The scalar constraint is thus compatible with the
semidiscrete BGN equations; its role becomes essential after time
discretization.

\begin{remark}[Compatibility of the initial polygon]\label{rem:compatibility}
Let $\ell_j=|\mathbf h_j|$, $\mathbf t_j=\mathbf h_j/\ell_j$, and
$\mathbf b_j=-(\mathbf h_j+\mathbf h_{j+1})^\perp/2$, with cyclic indices.
Testing \eqref{eq:semib} with the vector-valued nodal basis functions at
vertex $j$ gives
$\kappa_j\mathbf b_j=\mathbf t_j-\mathbf t_{j+1}$.
Taking its dot product with $\mathbf h_j+\mathbf h_{j+1}$ yields
\begin{equation}\label{eq:mesh-compatibility}
(\ell_j-\ell_{j+1})(1-\mathbf t_j\cdot\mathbf t_{j+1})=0.
\end{equation}
Thus adjacent edges with different unit tangents must have equal lengths.
When $\mathbf b_j\ne0$, the Euclidean least-squares approximation to the
nodal curvature equation is
\begin{equation}\label{eq:ls-curvature}
\kappa_j^{\mathrm{LS}}
=\frac{\mathbf b_j\cdot(\mathbf t_j-\mathbf t_{j+1})}{|\mathbf b_j|^2}.
\end{equation}
This value can be used to initialize the fully discrete curvature, but it
does not enforce \eqref{eq:semib} on an incompatible polygon. Accordingly,
the consistency analysis below assumes the existence of a smooth,
compatible semidiscrete reference solution.
\end{remark}

\subsection{BDF discretization with an energy-derivative constraint}

Following the BGN/BDF construction of \cite{jiang2024stable} and the
Lagrange-multiplier approach of \cite{garcke2025structure}, we first
apply the BDF operator to both the position and the energy.
Let $t_m=m\tau$, where $\tau=T/M>0$ and $M\in\mathbb N$.
The predicted polygon $\widetilde\Gamma^{m+1}$ is constructed before
assembling the equations and has the same nodal connectivity as the
computed polygons. Functions on different polygons are identified by
their coefficients on the common reference mesh. The normal and
arc-length derivatives in the following equations are evaluated on
$\widetilde\Gamma^{m+1}$. Given the required history, for $m\geq k-1$ seek
$(\mathbf X^{m+1},\kappa^{m+1},\lambda^{m+1})\in[V^h]^2\times V^h\times\mathbb R$
such that, for all $(\varphi^h,\boldsymbol\omega^h)\in V^h\times[V^h]^2$,
\begin{subequations}
	\begin{align}
		&\left(\frac{D_t^k \mathbf{X}^{m+1}}{\tau},
		\,\varphi^h \widetilde{\mathbf{n}}^{m+1}\right)^h_
		{\widetilde{\Gamma}^{m+1}}
		+ \left(\kappa^{m+1},\varphi^h\right)^h_
		{\widetilde{\Gamma}^{m+1}}
		- \lambda^{m+1}
		\left(\kappa^{m+1},\varphi^h\right)^h_
		{\widetilde{\Gamma}^{m+1}} = 0,
		\label{eq:discrete-evolution} \\
		&\left(\kappa^{m+1}\widetilde{\mathbf{n}}^{m+1},
		\,\boldsymbol{\omega}^h\right)^h_{\widetilde{\Gamma}^{m+1}}
		-\left(\partial_s\mathbf{X}^{m+1},
		\partial_s\boldsymbol{\omega}^h\right)_
		{\widetilde{\Gamma}^{m+1}}=0,
		\label{eq:discrete-geometric} \\
		&\frac{D_t^k W^{m+1}}{\tau}
		=-(\kappa^{m+1},\kappa^{m+1})^h_
		{\widetilde{\Gamma}^{m+1}}.
		\label{eq:energy dissipation D_t^k}
	\end{align}
\end{subequations}
The dimensionless BDF operators are
\[
\begin{aligned}
	D_t^1\mathbf X^{m+1}
	&= \mathbf X^{m+1}-\mathbf X^m,
	\\
	D_t^2\mathbf X^{m+1}
	&= \frac32\mathbf X^{m+1}
	-2\mathbf X^m
	+\frac12\mathbf X^{m-1},
	\\
	D_t^3\mathbf X^{m+1}
	&= \frac{11}{6}\mathbf X^{m+1}
	-3\mathbf X^m
	+\frac32\mathbf X^{m-1}
	-\frac13\mathbf X^{m-2},
	\\
	D_t^4\mathbf X^{m+1}
	&= \frac{25}{12}\mathbf X^{m+1}
	-4\mathbf X^m
	+3\mathbf X^{m-1}
	-\frac43\mathbf X^{m-2}
	+\frac14\mathbf X^{m-3}.
\end{aligned}
\]
The same coefficients define $D_t^kW^{m+1}$, with
$W^m:=W(\mathbf X^m):=L(\Gamma^m)$. Thus $D_t^k/\tau$
approximates the time derivative. We call these schemes BGN-LM-BDF\(k\),
$k=1,2,3,4$. For $k=1$, set $\widetilde\Gamma^{m+1}=\Gamma^m$;
for higher orders, use lower-order solves to predict the geometry.
Section~\ref{sec:staggered} specifies the predictors and start-up data.

\begin{remark}[Derivative-based energy updates]\label{rem:limitation}
Write $\Delta W^{m+1}=W^{m+1}-W^m$ and let $\mathcal D^{m+1}\geq0$
denote the endpoint dissipation. The BDF2 energy equation reads
\[
3\Delta W^{m+1}=\Delta W^m-2\tau\mathcal D^{m+1}.
\]
Consequently, stepwise decay follows by induction if $W^1\leq W^0$.
For BDF3 and BDF4, respectively,
\[
\begin{aligned}
11\Delta W^{m+1}
&=7\Delta W^m-2\Delta W^{m-1}-6\tau\mathcal D^{m+1},\\
25\Delta W^{m+1}
&=23\Delta W^m-13\Delta W^{m-1}+3\Delta W^{m-2}
-12\tau\mathcal D^{m+1}.
\end{aligned}
\]
Nonpositive history increments do not ensure a nonpositive right-hand
side. For BDF3, $W^0=2$, $W^1=W^2=1$, and $\mathcal D^3=0$
give $W^3=13/11>W^2$. For BDF4, $W^0=W^1=2$,
$W^2=W^3=1$, and $\mathcal D^4=0$ give $W^4=38/25>W^3$.
These examples concern only the scalar recurrences and illustrate the
limitation discussed in \cite{garcke2025structure}; they are not
counterexamples to a complete geometric scheme.
\end{remark}

The Crank--Nicolson formulation in
\cite{garcke2025structure} already uses an energy increment, with midpoint
dissipation. We instead couple an integrated energy balance to the
BDF approximation of the geometric evolution, using quadrature
consistent with the intended order.

\subsection{Energy-increment formulation and scalar analysis}

The proposed modification discretizes the scalar energy law in integral
form. Integrating
\eqref{eq:continuous-energy-dissipation} over
$[t_m,t_{m+1}]$ gives
\[
	W(t_{m+1})-W(t_m)
	+
	\int_{t_m}^{t_{m+1}}
	\int_{\Gamma(t)}\kappa^2\,ds\,dt=0.
\]
We retain the geometric equations
\eqref{eq:discrete-evolution}--\eqref{eq:discrete-geometric} and replace
\eqref{eq:energy dissipation D_t^k} by
\begin{equation}\label{eq:energy dissipation I_h^k}
	W^{m+1}-W^m+\tau\mathcal{Q}_k^{m+1}=0,
\end{equation}
where $\mathcal Q_k^{m+1}$ approximates the time-averaged dissipation rate.
For $k=1$, the right-endpoint rule gives
	\begin{equation}\label{I_h^k1}
		\mathcal{Q}_1^{m+1}
		=
		(\kappa^{m+1},\kappa^{m+1})^h_
		{\widetilde{\Gamma}^{m+1}}.
	\end{equation}

For $k=2$, we use the trapezoidal rule,
	\begin{equation}\label{I_h^k2}
		\mathcal{Q}_2^{m+1}
		=
		\frac12\left[
		(\kappa^{m+1},\kappa^{m+1})^h_{\widetilde{\Gamma}^{m+1}}
		+
		(\kappa^m,\kappa^m)^h_{\Gamma^m}
		\right].
	\end{equation}

For $k=3,4$, we use Simpson's rule,
	\begin{equation}\label{I_h^k34}
		\mathcal{Q}_k^{m+1}
		=
		\frac16\left[
		(\kappa^{m+1},\kappa^{m+1})^h_{\widetilde{\Gamma}^{m+1}}
		+4(\kappa^{m+\frac12},\kappa^{m+\frac12})^h_
		{\Gamma^{m+\frac12}}
		+(\kappa^m,\kappa^m)^h_{\Gamma^m}
		\right].
	\end{equation}
The resulting schemes are denoted by ES-BDF\(k\).
ES-BDF3 and ES-BDF4 share this quadrature rule but differ in their BDF
operators, predictors, and accuracy requirements for the intermediate data.
The endpoint dissipation is evaluated on the predicted geometry,
whereas the history and midpoint terms use their computed geometries.
We now examine consistency and energy dissipation for this scalar update.

For the consistency analysis, assume that \eqref{eq:semi} admits a
sufficiently smooth, compatible reference solution and write
\[
	\mathcal D^h(t)
	:=
	(\kappa^h(t),\kappa^h(t))_{\Gamma^h(t)}^h .
\]
The semidiscrete energy law \eqref{eq:semic} gives
\begin{equation}\label{eq:exact-averaged-energy-law}
	\frac{W^h(t_{m+1})-W^h(t_m)}{\tau}
	+
	\frac1{\tau}\int_{t_m}^{t_{m+1}}\mathcal D^h(t)\,dt
	=0 .
\end{equation}

Denote by $\mathcal D_{m,\theta}$ the computed dissipation at quadrature
node $t_{m+\theta}$ for the interval $[t_m,t_{m+1}]$. In the present
curve-shortening formulation,
\[
	\mathcal D_{m,1}
	=
	(\kappa^{m+1},\kappa^{m+1})^h_{\widetilde{\Gamma}^{m+1}},
	\qquad
	\mathcal D_{m,0}
	=
	(\kappa^{m},\kappa^{m})^h_{\Gamma^{m}},
\]
and, when Simpson's rule is used,
\[
	\mathcal D_{m,1/2}
	=
	(\kappa^{m+\frac12},\kappa^{m+\frac12})^h_
	{\Gamma^{m+\frac12}}.
\]

The required accuracy of these data is expressed by
\begin{equation}\label{eq:nodal-dissipation-consistency}
	\mathcal D_{m,\theta}
	=
	\mathcal D^h(t_{m+\theta})+O(\tau^k)
\end{equation}
at all nodes of the selected rule, uniformly in time for fixed $h$.
This is an assumption on both curvature and geometry, since the latter
determines the mass-lumped weights. In particular,
$\mathcal D_{m,0}$ need not equal $\mathcal D_{m-1,1}$: the same curvature
coefficients are evaluated on $\Gamma^m$ when used as history, but on
$\widetilde\Gamma^m$ when used as an endpoint.
For $k=3,4$, the assumption also concerns the midpoint states supplied by
the staggered construction in Section~\ref{sec:staggered}.

\begin{proposition}[Conditional scalar consistency]\label{prop:consistency}
Let a compatible semidiscrete solution exist on
$\mathcal I_m=[t_m,t_{m+1}]$, with $\mathcal D^h\in C^q(\mathcal I_m)$.
For the right-endpoint, trapezoidal,
and Simpson rules, respectively, set
\[
(q,C_q)=(1,1/2),\quad(2,1/12),\quad(4,1/2880).
\]
Define $e_m=\max_\theta|\mathcal D_{m,\theta}
-\mathcal D^h(t_{m+\theta})|$, with the maximum taken over the rule's nodes.
Then the defect in the energy-increment equation satisfies
\begin{equation}\label{eq:energy-increment-consistency}
\left|W^h(t_{m+1})-W^h(t_m)+\tau\mathcal Q_k^{m+1}\right|
\leq \tau e_m+C_q\tau^{q+1}
\|(\mathcal D^h)^{(q)}\|_{L^\infty(\mathcal I_m)}.
\end{equation}
Under \eqref{eq:nodal-dissipation-consistency}, this defect is
$O(\tau^{k+1})$ for $k=1,2,3,4$, or $O(\tau^k)$ after division by $\tau$,
with constants at fixed $h$.
\end{proposition}
\begin{proof}
By \eqref{eq:exact-averaged-energy-law}, the defect is the difference
between $\tau\mathcal Q_k^{m+1}$ and the exact dissipation integral.
Add and subtract the quadrature rule applied to the exact dissipation
values. If $\omega_\theta$ are the normalized
weights, the data error is bounded by
$\tau\sum_\theta\omega_\theta
|\mathcal D_{m,\theta}-\mathcal D^h(t_{m+\theta})|\leq\tau e_m$,
since $\omega_\theta\geq0$ and $\sum_\theta\omega_\theta=1$.
The standard quadrature remainder gives the second term in
\eqref{eq:energy-increment-consistency}. Finally, $q\geq k$ for each
selected rule, which gives the stated orders when $e_m=O(\tau^k)$.
\end{proof}
Proposition~\ref{prop:consistency} concerns the scalar constraint with
exact reference energies at the interval endpoints. A convergence result
for the coupled algorithm would also require control of the geometry,
curvature, initialization, predictors, and staggered histories; the
diffusion extensions additionally involve discrete-gradient errors.
The estimate is at fixed $h$, with no assertion of uniform constants
under spatial refinement or mesh degeneration.

\begin{proposition}[Stepwise energy dissipation]\label{prop:energy}
Let $\tau>0$, let all quadrature geometries be nondegenerate, and define
$\mathcal Q_k^{m+1}$ by \eqref{I_h^k1}--\eqref{I_h^k34}.
Every solution satisfying the scalar equation
\eqref{eq:energy dissipation I_h^k} exactly obeys
	
	\[
		W^{m+1}\le W^m.
	\]
For curve shortening flows, the inequality is strict if
\eqref{eq:discrete-geometric} also holds and the updated polygon has
at least one nonzero edge.
\end{proposition}

\begin{proof}
Each dissipation term is nonnegative, and each quadrature weight is
positive. Hence $\mathcal Q_k^{m+1}\geq0$, and
\eqref{eq:energy dissipation I_h^k} gives
	
	\[
		W^{m+1}-W^pm=-\tau\mathcal Q_k^{m+1}\le0.
	\]
Equality would force the endpoint dissipation to vanish. For curve shortening, positivity of
the mass-lumped weights on $\widetilde\Gamma^{m+1}$ then gives
$\kappa^{m+1}=0$. Taking
$\boldsymbol\omega^h=\mathbf X^{m+1}$ in
\eqref{eq:discrete-geometric} would then give
$\|\partial_s\mathbf X^{m+1}\|^2_{L^2(\widetilde\Gamma^{m+1})}=0$,
so every edge of the updated polygon would vanish, a contradiction.
\end{proof}

The proposition applies to every solved step, but does not establish
solvability for arbitrary $\tau$. In particular, the Simpson update
and $W^{m+1}>0$ require
\begin{equation}\label{eq:necessary-solvability}
\frac{\tau}{6}\bigl(4\mathcal D_{m,1/2}+\mathcal D_{m,0}\bigr)<W^m.
\end{equation}
Thus the prescribed history and midpoint data can preclude a
positive-energy solution; \eqref{eq:necessary-solvability} is only a
necessary condition. For successive solved steps on a single time grid,
the scalar equations telescope to
\[
W^M+\tau\sum_{m=m_0}^{M-1}\mathcal Q_k^{m+1}=W^{m_0}.
\]

This identity holds separately on the two grids used by the staggered
implementation. For inexact nonlinear solves, the corresponding residual
correction is discussed in Section~\ref{sec:newton}.
The principal role of the increment constraint is at BDF orders three
and four: BDF2 with a nonincreasing start, Crank--Nicolson, and the
trapezoidal increment update already admit energy-decay arguments.

Boole's rule has positive weights $(7,32,12,32,7)/90$, but a higher-order
extension would also require accurate intermediate states and a stable
multistep coupling. Not all closed Newton--Cotes rules have positive
weights. We focus here on design orders $k=2,3,4$, with ES-BDF1 as a
first-order reference. The remaining implementation task is to supply the midpoint geometry and
curvature required by Simpson's rule. Section~\ref{sec:staggered} describes
the interleaved time grids and the associated accuracy requirements.

\section{Staggered implementation}\label{sec:staggered}

We now describe the implementation of the ES-BDF\(k\) schemes for
curve-shortening flow. For ES-BDF3 and ES-BDF4, Simpson's rule requires
geometric data at the midpoint of each time interval. We obtain these data
by evolving an auxiliary sequence on the half-integer time grid.

Two ingredients follow the BGN/BDF\(k\) methodology of
\cite{jiang2024stable}. First, the usual start-up target for an order-$k$
multistep method is an $O(\tau^k)$ error in the position history. A first-order start-up
integrator with substep $\widetilde\tau$ over an interval of length
$O(\tau)$ has an accumulated error of order $O(\tau\widetilde\tau)$,
provided the local error and stability constants remain bounded. This
motivates $\widetilde\tau=O(\tau^{k-1})$. Choose
$\widetilde\tau=\tau^{k-1}/(2L)$ with an integer $L$ so that every half-grid
start-up time is reached exactly. This scaling does not by itself control
the algebraic curvature variable or verify
\eqref{eq:nodal-dissipation-consistency}. A first-order scheme may have
an $O(\widetilde\tau)$ curvature error even when its accumulated position
error is $O(\tau\widetilde\tau)$. The effects of these distinct errors
on the coupled start-up require separate analysis.

Second, each target geometry is predicted before the BGN equations
are assembled. We use BGN1 for ES-BDF2 and an ES-BDF2 predictor for
ES-BDF3. Here BGN1 means \eqref{eq:discrete-evolution} and
\eqref{eq:discrete-geometric} with $k=1$, $\lambda=0$, and no additional
scalar energy equation; it is distinct from ES-BDF1. Each nested
ES-BDF2 prediction uses its own BGN1 geometry.
ES-BDF2 uses only endpoint data and requires one time grid.
For the Simpson updates, Algorithm~\ref{alg:ES-BDF3} advances integer
and half-integer sequences, each with time step $\tau$.

\begin{algorithm}[htbp]
	\caption{ES-BDF3 staggered update}
	\label{alg:ES-BDF3}
	\begin{algorithmic}[1]
		\REQUIRE Initial data $\mathbf{X}^0$, time step $\tau$,
		start-up step $\widetilde{\tau}=\tau^2/(2L)=O(\tau^2)$,
        terminal time $T=M\tau$ with integer $M\geq3$.
		\ENSURE Integer-level states
		$(\mathbf{X}^m,\kappa^m)$ and half-integer-level states
		$(\mathbf{X}^{m+\frac12},\kappa^{m+\frac12})$.

		\STATE \textbf{Initialization.}
        \STATE Compute $\kappa^0$ by \eqref{eq:ls-curvature}.
        \STATE Compute $(\mathbf{X}^{1/2},\kappa^{1/2})$,
		$(\mathbf{X}^{1},\kappa^{1})$, and
		$(\mathbf{X}^{3/2},\kappa^{3/2})$ by successive BGN1
		start-up solves with step $\widetilde{\tau}$.
		\STATE Compute $(\mathbf{X}^{2},\kappa^{2})$ on the integer grid
		by ES-BDF2.
		\STATE Compute $(\mathbf{X}^{5/2},\kappa^{5/2})$ on the
		half-integer grid by ES-BDF2.

		\STATE \textbf{Time stepping.}
		\FOR{$m=2,\dots,M-1$}
		\STATE Construct the predicted geometry
		$\widetilde{\Gamma}^{m+1}$ by an ES-BDF2 predictor on the
		integer grid.
		\STATE Solve the ES-BDF3 system for
		$(\mathbf{X}^{m+1},\kappa^{m+1},\lambda^{m+1})$.
		\IF{$m<M-1$}
        \STATE Construct the predicted geometry
        $\widetilde{\Gamma}^{m+\frac32}$ by an ES-BDF2 predictor on
		the half-integer grid.
		\STATE Solve the ES-BDF3 system for
		$(\mathbf{X}^{m+\frac32},\kappa^{m+\frac32},
		\lambda^{m+\frac32})$.
        \ENDIF
        \ENDFOR
	\end{algorithmic}
\end{algorithm}

Algorithm~\ref{alg:ES-BDF3} advances the integer and half-integer sequences
alternately. The integer-level update at $t_{m+1}$ uses the previously
computed midpoint data at $t_{m+\frac12}$. The subsequent half-integer update
at $t_{m+\frac32}$ then uses the newly computed state at $t_{m+1}$ as its
midpoint data. The corresponding scalar energy equation on the half-integer
grid is
\begin{equation}
	W^{m+\frac32}-W^{m+\frac12}
	=
	-\tau\mathcal Q_3^{m+\frac32},
\end{equation}
where
\begin{equation}
	\mathcal Q_3^{m+\frac32}
	=
	\frac16\left[
	(\kappa^{m+\frac32},\kappa^{m+\frac32})^h_
	{\widetilde{\Gamma}^{m+\frac32}}
	+
	4(\kappa^{m+1},\kappa^{m+1})^h_{\Gamma^{m+1}}
	+
	(\kappa^{m+\frac12},\kappa^{m+\frac12})^h_
	{\Gamma^{m+\frac12}}
	\right].
\end{equation}

After initialization, both sequences use BDF3. This matches their
design orders but leaves \eqref{eq:nodal-dissipation-consistency} to be
verified for the coupled evolution. Proposition~\ref{prop:energy} applies
separately to the two sequences: it gives
$W^{m+1}\leq W^m$ and $W^{m+3/2}\leq W^{m+1/2}$, without ordering
$W^{m+1/2}$ and $W^m$.

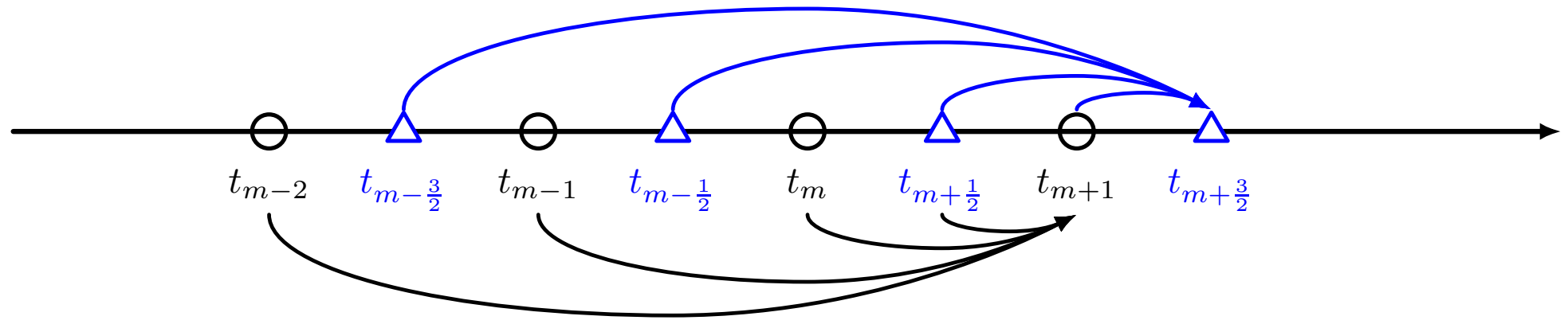
\begin{figure}[htbp]
	\centering
	\resizebox{0.99\textwidth}{!}{%
		\begin{tikzpicture}[
			>={Latex[length=2mm,width=1.6mm]},
			line width=1.05pt,
			line cap=round,
			line join=round,
			every node/.style={font=\normalsize, inner sep=1.5pt}
			]

			\def\dx{1.25}
			\def\circler{4.5pt}
			\def\triSize{10pt}
			\def\blue{blue}
			\def\blackdrop{22pt}
			\def\triTop{5.8pt}
			\def\labelsep{8pt}
			\def\kappa{0.55228475}

			\draw[line width=1.35pt, ->] (-5.9*\dx,0) -- (5.6*\dx,0);

			\coordinate (tm2)  at (-4*\dx,0);
			\coordinate (tm3h) at (-3*\dx,0);
			\coordinate (tm1)  at (-2*\dx,0);
			\coordinate (tm1h) at (-1*\dx,0);
			\coordinate (tm)   at ( 0,0);
			\coordinate (tp1h) at ( 1*\dx,0);
			\coordinate (tp1)  at ( 2*\dx,0);
			\coordinate (tp3h) at ( 3*\dx,0);

			\draw[line width=1.1pt] (tm2) circle (\circler);
			\draw[line width=1.1pt] (tm1) circle (\circler);
			\draw[line width=1.1pt] (tm)  circle (\circler);
			\draw[line width=1.1pt] (tp1) circle (\circler);

			\node[below=\labelsep of tm2] {$t_{m-2}$};
			\node[below=\labelsep of tm1] {$t_{m-1}$};
			\node[below=\labelsep of tm]  {$t_m$};
			\node[below=\labelsep of tp1] {$t_{m+1}$};

			\node[draw=\blue, regular polygon, regular polygon sides=3,
			minimum size=\triSize, fill=white, inner sep=0pt] at (tm3h) {};
			\node[draw=\blue, regular polygon, regular polygon sides=3,
			minimum size=\triSize, fill=white, inner sep=0pt] at (tm1h) {};
			\node[draw=\blue, regular polygon, regular polygon sides=3,
			minimum size=\triSize, fill=white, inner sep=0pt] at (tp1h) {};
			\node[draw=\blue, regular polygon, regular polygon sides=3,
			minimum size=\triSize, fill=white, inner sep=0pt] at (tp3h) {};

			\node[below=\labelsep of tm3h, \blue]
			{$t_{m-\frac{3}{2}}$};
			\node[below=\labelsep of tm1h, \blue]
			{$t_{m-\frac{1}{2}}$};
			\node[below=\labelsep of tp1h, \blue]
			{$t_{m+\frac{1}{2}}$};
			\node[below=\labelsep of tp3h, \blue]
			{$t_{m+\frac{3}{2}}$};

			\coordinate (tm3hTop) at ($(tm3h)+(0,\triTop)$);
			\coordinate (tm1hTop) at ($(tm1h)+(0,\triTop)$);
			\coordinate (tp1hTop) at ($(tp1h)+(0,\triTop)$);
			\coordinate (tp3hTop) at ($(tp3h)+(0,\triTop)$);
			\coordinate (tp1BlueTop) at ($(tp1)+(0,\triTop)$);

			\coordinate (tm2Bot) at ($(tm2)+(0,-\blackdrop)$);
			\coordinate (tm1Bot) at ($(tm1)+(0,-\blackdrop)$);
			\coordinate (tmBot) at ($(tm)+(0,-\blackdrop)$);
			\coordinate (tp1Bot) at ($(tp1)+(0,-\blackdrop)$);
			\coordinate (tp1hBot) at ($(tp1h)+(0,-\blackdrop)$);

			\newcommand{\UpperEllipseArrow}[3]{%
				\draw[#3, ->, line width=1pt]
				let \p1 = (#1), \p2 = (#2),
				\n1 = {(\x2-\x1)/2},
				\n2 = {\n1/4},
				\n3 = {\kappa*\n1},
				\n4 = {\kappa*\n2}
				in
				(#1)
				.. controls ($ (#1) + (0,\n4) $)
				and ($ (#1) + (\n1-\n3,\n2) $)
				.. ($ ({(\x1+\x2)/2},{\y1+\n2}) $)
				.. controls
				($ ({(\x1+\x2)/2+\n3},{\y1+\n2}) $)
				and ($ (#2) + (-4pt,2pt) $)
				.. ($ (#2) + (-4pt,2pt) $)
				-- (#2);
			}

			\newcommand{\LowerEllipseArrow}[3]{%
				\draw[#3, ->, line width=1pt]
				let \p1 = (#1), \p2 = (#2),
				\n1 = {(\x2-\x1)/2},
				\n2 = {\n1/4},
				\n3 = {\kappa*\n1},
				\n4 = {\kappa*\n2}
				in
				(#1)
				.. controls ($ (#1) + (0,-\n4) $)
				and ($ (#1) + (\n1-\n3,-\n2) $)
				.. ($ ({(\x1+\x2)/2},{\y1-\n2}) $)
				.. controls
				($ ({(\x1+\x2)/2+\n3},{\y1-\n2}) $)
				and ($ (#2) + (-4pt,-2pt) $)
				.. ($ (#2) + (-4pt,-2pt) $)
				-- (#2);
			}

			\UpperEllipseArrow{tm3hTop}{tp3hTop}{\blue}
			\UpperEllipseArrow{tm1hTop}{tp3hTop}{\blue}
			\UpperEllipseArrow{tp1hTop}{tp3hTop}{\blue}
			\UpperEllipseArrow{tp1BlueTop}{tp3hTop}{\blue}

			\LowerEllipseArrow{tm2Bot}{tp1Bot}{black}
			\LowerEllipseArrow{tm1Bot}{tp1Bot}{black}
			\LowerEllipseArrow{tmBot}{tp1Bot}{black}
			\LowerEllipseArrow{tp1hBot}{tp1Bot}{black}

		\end{tikzpicture}%
	}
	\caption{	Dependency structure of one staggered ES-BDF3 update.
	Black arrows indicate data used for the integer-level update at
	$t_{m+1}$, and blue arrows indicate data used for the subsequent
	half-integer-level update at $t_{m+\frac32}$. The arrows include both
	the BDF3 history and the midpoint data required by Simpson's rule.
	}
	\label{fig:ESBDF3-timeline}
\end{figure}

The dependency structure of one staggered update is illustrated in
Fig.~\ref{fig:ESBDF3-timeline}.

For ES-BDF4, take $M\geq4$ and retain four history values on each grid. The integer sequence is initialized
through $t_3$ and the half-integer sequence through $t_{7/2}$ using
lower-order solves refined sufficiently to target $O(\tau^4)$ start-up
errors. An unrefined ES-BDF2 macro step generally contributes an
$O(\tau^3)$ local error and is not by itself sufficient for this purpose.
The ES-BDF3 predictor uses only the already available history and midpoint
data and is discarded after construction of the predicted geometry. The histories are listed in Table~\ref{tab:ESBDFk-implementation}.
As for ES-BDF3, the integer solve uses the available half-integer midpoint,
and the following half-integer solve uses the newly computed integer state.
Each ES-BDF3 or ES-BDF4 macro step normally
requires two nonlinear solves, together with its predictor solves. Newton
iteration counts alone therefore do not measure total computational cost.

\begin{table}[htbp]
	\centering
	\caption{Staggered implementation data for ES-BDF3 and ES-BDF4.}
	\label{tab:ESBDFk-implementation}
	\small
	\setlength{\tabcolsep}{5pt}
	\renewcommand{\arraystretch}{1.15}
	\begin{tabular}{c|c|c}
		\hline
		& ES-BDF3 & ES-BDF4
		\\
		\hline
		Integer-level history
		& \(m,m-1,m-2\)
		& \(m,m-1,m-2,m-3\)
		\\
		Half-integer-level history
		& \(m+\frac12,m-\frac12,m-\frac32\)
		& \(m+\frac12,m-\frac12,m-\frac32,m-\frac52\)
		\\
		Predictor
		& ES-BDF2
		& ES-BDF3
		\\
		Integer startup
		& through \(t_2\)
		& through \(t_3\)
		\\
		Half-integer startup
		& through \(t_{5/2}\)
		& through \(t_{7/2}\)
		\\
		\hline
	\end{tabular}
\end{table}

\section{Newton solution of the ES-BDF3 system}\label{sec:newton}

At the target time $t_{m+1}$, Newton's method solves for
$(\mathbf X^{m+1},\kappa^{m+1},\lambda^{m+1})
\in[V^h]^2\times V^h\times\mathbb R$.
The predicted geometry, its normal and metric, and the midpoint data are
fixed throughout the iteration. Initial curvature is obtained from
\eqref{eq:ls-curvature}, following
\cite[Remark~2.1]{jiang2024second}; this requires
$\mathbf b_j^0\ne0$, which is stronger than nonzero individual edge lengths.
The position and curvature histories are supplied by
Section~\ref{sec:staggered}.

For a polygonal curve $\Gamma=\mathbf X(\mathbb I)$ and a perturbation
$Y(\varepsilon)=\mathbf X+\varepsilon\mathbf X^\delta$, the first variation
of the discrete perimeter is
\begin{equation}
	\delta W(\mathbf{X})(\mathbf{X}^\delta)
	=
	\left.
	\frac{d}{d\varepsilon}W(Y(\varepsilon))
	\right|_{\varepsilon=0}
	=
	(\partial_s\mathbf X,\partial_s\mathbf X^\delta)_\Gamma .
	\label{eq:varW}
\end{equation}
Equivalently, $\delta W(\mathbf X)(\mathbf X^\delta)
=\sum_j (\mathbf h_j/|\mathbf h_j|)\cdot
(\mathbf X^\delta(\rho_j)-\mathbf X^\delta(\rho_{j-1}))$.
This derivative is evaluated on the current iterate $\Gamma^{m+1,i}$,
whereas the metric in the endpoint dissipation remains fixed on
$\widetilde\Gamma^{m+1}$.

A natural initial guess is
\[
	\mathbf{X}^{m+1,0}=\widetilde{\mathbf{X}}^{m+1},
	\qquad
	\kappa^{m+1,0}=\kappa^{m},
	\qquad
	\lambda^{m+1,0}=0.
\]

Given
$(\mathbf{X}^{m+1,i},\kappa^{m+1,i},\lambda^{m+1,i})$, the Newton correction
$(\mathbf{X}^\delta,\kappa^\delta,\lambda^\delta)
\in[V^h]^2\times V^h\times\mathbb R$
is determined by the following linearized system. Let
$R_1^i(\varphi^h)$ and $R_2^i(\boldsymbol\omega^h)$ be the left-hand sides
of \eqref{eq:discrete-evolution} and \eqref{eq:discrete-geometric}, evaluated
at the current iterate with $k=3$, and let
$R_E^i=W(\mathbf X^{m+1,i})-W^m+\tau\mathcal Q_3^{m+1,i}$.
Then
\begin{align}
&\left(\frac{11}{6\tau}\mathbf X^\delta,
 \varphi^h\widetilde{\mathbf n}^{m+1}\right)^h_{\widetilde\Gamma^{m+1}}
 +(1-\lambda^{m+1,i})(\kappa^\delta,\varphi^h)^h_{\widetilde\Gamma^{m+1}}
 \nonumber\\
&\qquad-\lambda^\delta(\kappa^{m+1,i},\varphi^h)^h_{\widetilde\Gamma^{m+1}}
 =-R_1^i(\varphi^h),\label{eq:Newton-1}\\
&(\kappa^\delta\widetilde{\mathbf n}^{m+1},\boldsymbol\omega^h)^h_{\widetilde\Gamma^{m+1}}
 -(\partial_s\mathbf X^\delta,\partial_s\boldsymbol\omega^h)_{\widetilde\Gamma^{m+1}}
 =-R_2^i(\boldsymbol\omega^h),\label{eq:Newton-2}\\
&(\partial_s\mathbf X^{m+1,i},\partial_s\mathbf X^\delta)_{\Gamma^{m+1,i}}
 +\frac\tau3(\kappa^{m+1,i},\kappa^\delta)^h_{\widetilde\Gamma^{m+1}}
 =-R_E^i.\label{eq:Newton-3}
\end{align}
The factor $\tau/3$ in \eqref{eq:Newton-3} is twice the endpoint Simpson
weight $\tau/6$; no geometry derivative enters this quadratic term because
the predicted metric is fixed.

\begin{proposition}[Regularity of the CSF Newton system]\label{prop:newton}
At an iterate with a nondegenerate current polygon and $\lambda^{m+1,i}<1$,
suppose that the predicted polygon
has positive edge lengths and its nodal vectors
$\widetilde{\mathbf b}_j
=-(\widetilde{\mathbf h}_j+\widetilde{\mathbf h}_{j+1})^\perp/2$
span $\mathbb R^2$. Let $K_i$ be the coefficient matrix of
\eqref{eq:Newton-1}--\eqref{eq:Newton-2} for
$z=(\mathbf X^\delta,\kappa^\delta)$, with $\lambda^\delta=0$.
Then $K_i$ is invertible. Writing the full Jacobian as
\begin{equation}\label{eq:newton-schur}
J_i=\begin{pmatrix}K_i&b_i\\c_i^\top&0\end{pmatrix},
\qquad s_i=-c_i^\top K_i^{-1}b_i,
\end{equation}
where $b_i$ is the multiplier column and $c_i^\top$ is the energy row,
$J_i$ is invertible if and only if $s_i\ne0$.
\end{proposition}
\begin{proof}
For a homogeneous solution of the first two equations, test with
$\varphi^h=\kappa^\delta$ and $\boldsymbol\omega^h=\mathbf X^\delta$.
Eliminating their mixed term gives
\[
\frac{11}{6\tau}\|\partial_s\mathbf X^\delta\|^2_{L^2(\widetilde\Gamma^{m+1})}
+(1-\lambda^{m+1,i})\|\kappa^\delta\|_{h,\widetilde\Gamma^{m+1}}^2=0,
\]
where $\|v\|_{h,\Gamma}^2=(v,v)^h_\Gamma$. Hence $\kappa^\delta=0$
and $\mathbf X^\delta$ is constant. The first equation and the spanning
condition force this constant to vanish. The square matrix $K_i$ thus
has trivial kernel. Block elimination then gives the scalar Schur
complement $s_i$ in \eqref{eq:newton-schur}, proving the second assertion.
\end{proof}
At a nondegenerate root satisfying these conditions, the inverse function
theorem gives local uniqueness of that root and standard Newton theory
gives local quadratic convergence for sufficiently close initial guesses.
This conditional statement proves neither existence of a root nor that
$s_i$ stays away from zero along the computed evolution. The condition
$\lambda^{m+1,i}<1$ is sufficient for the block argument, not asserted for
every iterate. These conclusions concern CSF, not the anisotropic matrix.

After solving \eqref{eq:Newton-1}--\eqref{eq:Newton-3}, the iterate is
updated by
\[
	\mathbf{X}^{m+1,i+1}
	=
	\mathbf{X}^{m+1,i}+\mathbf{X}^\delta,
	\qquad
	\kappa^{m+1,i+1}
	=
	\kappa^{m+1,i}+\kappa^\delta,
	\qquad
	\lambda^{m+1,i+1}
	=
	\lambda^{m+1,i}+\lambda^\delta.
\]

The increment-based stopping criterion used to describe the iteration is
\begin{equation}
	\label{eq:newton-stopping}
	\max\left\{
		\|\mathbf X^\delta\|_{L^\infty},
		\|\kappa^\delta\|_{L^\infty},
		|\lambda^\delta|
	\right\}
	\leq \mathrm{tol}.
\end{equation}
A small correction alone does not certify a small nonlinear residual.
In particular, if the accepted state has scalar energy residual
$r_E^{m+1}:=W^{m+1}-W^m+\tau\mathcal Q_k^{m+1}$, then
\begin{equation}\label{eq:energy-residual}
W^{m+1}-W^m=-\tau\mathcal Q_k^{m+1}+r_E^{m+1}.
\end{equation}
Consequently, a sufficient acceptance condition for dissipation is
$r_E^{m+1}\leq\tau\mathcal Q_k^{m+1}$. For example,
$|r_E^{m+1}|\leq\eta\tau\mathcal Q_k^{m+1}$ with $0\leq\eta<1$
gives $W^{m+1}-W^m\leq-(1-\eta)\tau\mathcal Q_k^{m+1}$.
When $\mathcal Q_k^{m+1}=0$, a positive absolute tolerance alone cannot
guarantee nonincrease. More generally,
\[
\left|W^M-W^{m_0}+\tau\sum_{m=m_0}^{M-1}\mathcal Q_k^{m+1}\right|
\leq\sum_{m=m_0}^{M-1}|r_E^{m+1}|.
\]
Thus a per-step energy residual $O(\tau^{k+1})$ contributes at most
$O(\tau^k)$ to the accumulated balance over a fixed time interval.
This controls the balance error, not the geometric solution error.
These estimates motivate monitoring the energy residual alongside
the correction size, with residual scaling, tolerances, and iteration
limits specified for each implementation. Control of the geometric
solution error additionally requires bounds on the other residuals.


\section{Extensions to other geometric flows}\label{sec:extensions}

The construction extends by changing the geometric constitutive relation
and the dissipation functional. To display the common structure, for a
nonnegative functional $d_\Gamma(u)$ define
\begin{equation}\label{eq:flow-quadrature}
\mathscr Q_m[d,u]
=\frac16\left[
d_{\widetilde\Gamma^{m+1}}(u^{m+1})
+4d_{\Gamma^{m+1/2}}(u^{m+1/2})
+d_{\Gamma^m}(u^m)\right].
\end{equation}
The endpoint uses the predicted metric, whereas the history and midpoint
use their computed geometries. The following ES-BDF3 systems share this
quadrature and the staggered construction of
Section~\ref{sec:staggered}. Nondegenerate evaluation geometries are assumed.

\subsection{Curve diffusion}

We first consider curve diffusion, the $H^{-1}$-gradient flow of the
perimeter $W(t)=L(\Gamma(t))$. Its normal velocity is
\[
	V=\partial_{ss}\kappa,
\]
where $\kappa$ denotes the curvature.
The corresponding Lagrange-multiplier formulation is
\begin{equation}\label{eq:CDF-continuous}
	\begin{aligned}
		\partial_t \mathbf{X}\cdot\mathbf{n}
		&= \partial_{ss}\kappa+\lambda(t)\kappa, \\
		\kappa\,\mathbf{n}
		&= -\partial_{ss}\mathbf{X}, \\
		\frac{dW}{dt}
		&= -\int_{\Gamma(t)} |\partial_s\kappa|^2\,ds.
	\end{aligned}
\end{equation}
Using the perimeter variation
$dW/dt=(\kappa,V)_{\Gamma(t)}$ and the first equation in
\eqref{eq:CDF-continuous} gives
\[
	\lambda(t)
	\int_{\Gamma(t)}\kappa^2\,ds=0.
\]
Hence $\lambda(t)=0$ whenever $\kappa\not\equiv0$, and the augmented
formulation is equivalent to the classical curve-diffusion equation at the
continuous level.

Define $d_{\mathrm C,\Gamma}(u)=\|\partial_su\|_{L^2(\Gamma)}^2$.
Given $\widetilde\Gamma^{m+1}$, the ES-BDF3 scheme seeks
\[
	(\mathbf X^{m+1},\kappa^{m+1},\lambda^{m+1})
	\in [V^h]^2\times V^h\times\mathbb R
\]
such that, for all
$(\varphi^h,\boldsymbol\omega^h)\in V^h\times[V^h]^2$,
\begin{subequations}\label{eq:ES-BDF3-CDF}
	\begin{align}
		\left(
		\frac{D_t^3\mathbf X^{m+1}}{\tau},
		\varphi^h \widetilde{\mathbf n}^{m+1}
		\right)^h_{\widetilde{\Gamma}^{m+1}}
		&+
		\left(
		\partial_s\kappa^{m+1},
		\partial_s\varphi^h
		\right)_{\widetilde{\Gamma}^{m+1}}
		\nonumber\\
		&-
		\lambda^{m+1}
		\left(
		\kappa^{m+1},\varphi^h
		\right)^h_{\widetilde{\Gamma}^{m+1}}
		=0,
		\\
		\left(
		\kappa^{m+1}\widetilde{\mathbf n}^{m+1},
		\boldsymbol{\omega}^h
		\right)^h_{\widetilde{\Gamma}^{m+1}}
		&-
		\left(
		\partial_s\mathbf X^{m+1},
		\partial_s\boldsymbol{\omega}^h
		\right)_{\widetilde{\Gamma}^{m+1}}
		=0,
		\\
		W^{m+1}-W^m
        &+\tau\mathscr Q_m[d_{\mathrm C},\kappa]=0.
	\end{align}
\end{subequations}

Since $\mathscr Q_m[d_{\mathrm C},\kappa]\geq0$, the last equation
gives $W^{m+1}\leq W^m$. Testing the first equation with $\varphi^h=1$ yields
\[
\left(\frac{D_t^3\mathbf X^{m+1}}{\tau},
\widetilde{\mathbf n}^{m+1}\right)^h_{\widetilde\Gamma^{m+1}}
=\lambda^{m+1}(\kappa^{m+1},1)^h_{\widetilde\Gamma^{m+1}},
\]
which is not an exact discrete enclosed-area conservation law. An additional
area constraint, as in \cite{garcke2025structure}, is needed if exact
conservation is required. 

\subsection{Anisotropic curve diffusion}

We next consider anisotropic curve diffusion driven by the anisotropic
perimeter
\[
	W_\gamma(\Gamma)=\int_\Gamma \gamma(\theta)\,ds,
\]
where $\theta$ is the angle of the outward normal measured counterclockwise
from the positive first coordinate axis. We assume a positive, $2\pi$-periodic
energy density $\gamma\in C^2$. The unregularized smooth-flow setting below
uses positive surface stiffness $\gamma+\gamma''>0$.
The normal velocity and anisotropic chemical potential are
\[
	V=\partial_{ss}\mu,
	\qquad
	\mu=(\gamma(\theta)+\gamma''(\theta))\kappa.
\]
For $\gamma\equiv1$, one has $\mu=\kappa$, and the formulation reduces to
isotropic curve diffusion. With the sign convention used here, the first
variation satisfies
\[
	\frac{d}{dt}W_\gamma(\Gamma(t))
	=
	(\mu,V)_{\Gamma(t)}.
\]

For a polygonal curve \(\Gamma^m\), we define the corresponding discrete
anisotropic energy by
\[
	W_\gamma^m
	:=
	W_\gamma(\Gamma^m)
	=
	\sum_{j=1}^{N}
	\gamma(\theta_j^m)\,|\mathbf h_j^m|,
\]
where \(\theta_j^m\) denotes the orientation angle of the outward unit normal
on the \(j\)th edge of \(\Gamma^m\).

For Newton iteration, the energy row must differentiate this anisotropic
energy. At the current polygon,
\begin{equation}\label{eq:anisotropic-variation}
\delta W_\gamma(\mathbf X)[\mathbf Y]
=\sum_j\bigl(\gamma(\theta_j)\mathbf t_j
-\gamma'(\theta_j)\mathbf n_j\bigr)\cdot
(\mathbf Y_j-\mathbf Y_{j-1}).
\end{equation}
Indeed, $\delta\ell_j=\mathbf t_j\cdot\delta\mathbf h_j$ and
$\delta\theta_j=-\mathbf n_j\cdot\delta\mathbf h_j/\ell_j$.
The energy derivative includes both edge-length and angle variations;
the coefficients in the predicted-geometry equations remain fixed.

Following \cite{li2021energy}, we use the conservative representation
\[
	\mu\,\mathbf{n}
	=
	-\partial_s\!\bigl(G(\theta)\,\partial_s\mathbf{X}\bigr),
	\qquad
	G(\theta)=
	\begin{pmatrix}
		\gamma(\theta) & -\gamma'(\theta)\\
		\gamma'(\theta) & \gamma(\theta)
	\end{pmatrix}.
\]
The augmented continuous formulation is then
\begin{equation}\label{eq:ACDF-continuous}
	\begin{aligned}
		\partial_t \mathbf{X}\cdot\mathbf{n}
		&= \partial_{ss}\mu+\lambda(t)\mu, \\
		\mu\,\mathbf{n}
		&= -\partial_s\bigl(G(\theta)\partial_s\mathbf{X}\bigr), \\
		\frac{dW_\gamma}{dt}
		&= -\int_{\Gamma(t)} |\partial_s\mu|^2\,ds.
	\end{aligned}
\end{equation}
Combining the first-variation identity with
\eqref{eq:ACDF-continuous} gives
\[
	\lambda(t)\int_{\Gamma(t)}\mu^2\,ds=0,
\]
and therefore $\lambda(t)=0$ whenever $\mu\not\equiv0$.

Using $d_{\mathrm C,\Gamma}(\mu)=\|\partial_s\mu\|_{L^2(\Gamma)}^2$
and the predicted geometry, the ES-BDF3 equations are, for every
$(\psi^h,\boldsymbol w^h)\in V^h\times[V^h]^2$,
\begin{subequations}\label{eq:ES-BDF3-ACDF}
	\begin{align}
		\left(
		\frac{D_t^3\mathbf X^{m+1}}{\tau},
		\widetilde{\mathbf n}^{m+1}\psi^h
		\right)^h_{\widetilde{\Gamma}^{m+1}}
		&+
		\left(
		\partial_s\mu^{m+1},
		\partial_s\psi^h
		\right)_{\widetilde{\Gamma}^{m+1}}
		\nonumber\\
		&-
		\lambda^{m+1}
		\left(
		\mu^{m+1},\psi^h
		\right)^h_{\widetilde{\Gamma}^{m+1}}
		=0,
		\\
		\left(
		\mu^{m+1}\widetilde{\mathbf n}^{m+1},
		\boldsymbol w^h
		\right)^h_{\widetilde{\Gamma}^{m+1}}
		&-
		\left(
		G(\widetilde{\theta}^{m+1})
		\partial_s\mathbf X^{m+1},
		\partial_s\boldsymbol w^h
		\right)_{\widetilde{\Gamma}^{m+1}}
		=0,
		\label{aniSDF b}
		\\
		W_\gamma^{m+1}-W_\gamma^m
        &+\tau\mathscr Q_m[d_{\mathrm C},\mu]=0.
	\end{align}
\end{subequations}

Again, the positivity of the Simpson weights gives
$W_\gamma^{m+1}\leq W_\gamma^m$ for any solution of
\eqref{eq:ES-BDF3-ACDF}.

\begin{remark}[Strong anisotropic cases]
The constitutive relation can also be expressed using the symmetrized
matrix of \cite{Bao2023Symmetrized,jiang2025predictor}. To distinguish the
stabilization parameter from the BDF order, we write this matrix as
$\mathbf Z_{\sigma}(\widetilde{\mathbf n}^{m+1})$ in place of
$G(\widetilde\theta^{m+1})$ in \eqref{aniSDF b}. The scalar energy-increment
argument still applies to any solution with nonnegative dissipation.
If $\gamma+\gamma''$ becomes negative in some orientations, however, algebraic energy decay
and matrix symmetrization alone do not establish well-posedness of the
unregularized continuum flow or convergence to it. Such a regime requires
a specified regularization or solution concept and is not discussed here.
\end{remark}

\paragraph{Initialization} We extend the least-squares construction of \cite[Remark~2.1]{jiang2024second} to the anisotropic
chemical potential $\mu$. For the initial polygon $\Gamma^0$, the resulting normal equations are
\[
	(N_0^\top N_0)\mu^0=N_0^\top A_0\mathbf X^0,
\]
where
\[
	N_0
	=
	\binom{
	(\varphi_i,(\mathbf n^0)^{[1]}\varphi_j)^h_{\Gamma^0}
	}{
	(\varphi_i,(\mathbf n^0)^{[2]}\varphi_j)^h_{\Gamma^0}
	},
	\qquad
	\mathbf X^0
	=
	\binom{\mathbf x^0}{\mathbf y^0},
\]
and
\[
A_0=\begin{pmatrix}A_{\gamma}&-A_{\gamma'}\\
                         A_{\gamma'}&A_{\gamma}\end{pmatrix},
\qquad
(A_\eta)_{ij}=(\eta(\theta^0)\partial_s\varphi_j,
                         \partial_s\varphi_i)_{\Gamma^0}.
\]
Here $\{\varphi_i\}_{i=1}^N$ is the standard Lagrange basis on $\mathbb I$,
$\mathbf a^{[1]}$ and $\mathbf a^{[2]}$ denote the two components of
$\mathbf a\in\mathbb R^2$, and
$x_j^0=(\mathbf X^0)^{[1]}(\rho_j)$ and
$y_j^0=(\mathbf X^0)^{[2]}(\rho_j)$.
Equivalently, $\mu^0$ minimizes
$\|N_0\mu-A_0\mathbf X^0\|_2^2$. Uniqueness requires $N_0$ to have
full column rank, or equivalently nonzero mass-weighted nodal normals.
If the system is ill conditioned, a QR or singular-value decomposition
avoids explicitly forming the normal equations. This supplies the initial chemical potential.

\subsection{Surface evolution in three dimensions}

We next consider closed surfaces in $\mathbb R^3$. The curvature relation is
expressed in terms of the mean curvature $\mathcal H$, and surface gradients
replace arc-length derivatives. We first introduce the discrete surface
notation used below.

At time $t_m=m\tau$, the surface $\Gamma(t_m)$ is approximated by a polyhedral
surface
\[
	\Gamma^m=\bigcup_{j=1}^{J}\sigma_j^m,\qquad
	\sigma_j^m
	=
	\triangle\{\mathbf{q}_{j1}^m,\mathbf{q}_{j2}^m,\mathbf{q}_{j3}^m\},
\]
consisting of nondegenerate triangles without hanging vertices, with
the vertices of each triangle ordered consistently so that their cross
product points outward. The continuous
piecewise linear finite element space on $\Gamma^m$ is
\[
	W_h^m
	:=
	\left\{
	u\in C(\Gamma^m):
	u|_{\sigma_j^m}\ \text{is affine for all }j=1,\dots,J
	\right\}
	\subset H^1(\Gamma^m).
\]

For scalar/vector piecewise continuous functions $u,v$ defined on $\Gamma^m$,
we denote
\[
	(u,v)_{\Gamma^m}
	:=
	\int_{\Gamma^m}u\cdot v\,dA,
	\qquad
	(u,v)_{\Gamma^m}^h
	:=
	\frac13\sum_{j=1}^{J}|\sigma_j^m|
	\sum_{k=1}^{3}
	(u\cdot v)\bigl((\mathbf{q}_{jk}^m)^-\bigr),
\]
where $(\mathbf{q}_{jk}^m)^-$ denotes the limit taken from within the triangle
$\sigma_j^m$.

The outer unit normal on each triangle is given by
\[
	\mathbf{n}^m|_{\sigma_j^m}
	=
	\mathbf{n}_j^m
	:=
	\frac{
	(\mathbf{q}_{j2}^m-\mathbf{q}_{j1}^m)
	\times
	(\mathbf{q}_{j3}^m-\mathbf{q}_{j1}^m)
	}{
	\bigl|
	(\mathbf{q}_{j2}^m-\mathbf{q}_{j1}^m)
	\times
	(\mathbf{q}_{j3}^m-\mathbf{q}_{j1}^m)
	\bigr|
	},
\]
and the discrete surface gradient $\nabla_{\Gamma^m} f$ of $f\in W_h^m$ is
constant on each triangle and defined by
\[
	(\nabla_{\Gamma^m} f)|_{\sigma_j^m}
	=
	\frac{
    \mathbf n_j^m\times\bigl(
    f(\mathbf{q}_{j1}^m)(\mathbf{q}_{j3}^m-\mathbf{q}_{j2}^m)
	+
	f(\mathbf{q}_{j2}^m)(\mathbf{q}_{j1}^m-\mathbf{q}_{j3}^m)
	+
	f(\mathbf{q}_{j3}^m)(\mathbf{q}_{j2}^m-\mathbf{q}_{j1}^m)
	\bigr)
    }{
    2|\sigma_j^m|
	}.
\]
The discrete surface area of $\Gamma^m$ is denoted by
\[
	W^m:=S^m:=\int_{\Gamma^m}1\,dA.
\]
For a nondegenerate current polyhedral surface $\Gamma$ with vertex vector $\mathbf X$,
its area variation in direction $\mathbf Y$ is
\begin{equation}\label{eq:surface-area-variation}
\delta S(\mathbf X)[\mathbf Y]
=\int_\Gamma\nabla_\Gamma\mathrm{Id}:\nabla_\Gamma\mathbf Y\,dA.
\end{equation}
Here $\mathbf Y$ is identified with its field on the current surface.
This is the surface analogue of \eqref{eq:varW} for the Newton energy row;
the gradients in the dissipation term remain on the frozen predictor.

\paragraph{Mean curvature flow and surface diffusion}

Let $\mathcal H$ denote the sum of the two principal curvatures (rather
than their average), with the sign convention
\[
	\mathcal{H}\,\mathbf{n}
	=
	-\Delta_\Gamma \mathrm{Id}.
\]
For mean curvature flow,
$V=-\mathcal H$, and the surface area
$W(t)=S(\Gamma(t))$ satisfies
\[
	\frac{dW}{dt}
	=
	-\int_{\Gamma(t)}|\mathcal{H}|^2\,dA.
\]
The augmented velocity is $V=-(1-\lambda)\mathcal H$, with the same
curvature relation and area-dissipation law. The first variation
$dW/dt=(\mathcal H,V)_\Gamma$ then gives
$\lambda\|\mathcal H\|_{L^2(\Gamma)}^2=0$, so the augmentation reduces to
mean curvature flow whenever $\mathcal H\not\equiv0$.

For mean curvature flow, set
$d_{\mathrm M,\Gamma}(u)=(u,u)^h_\Gamma$ in
\eqref{eq:flow-quadrature}. All surfaces share a fixed mesh connectivity. The map
$\mathbf X^{m+1}$ and historical position fields are represented on the
predicted surface by nodal identification; scalar test functions belong to
its continuous piecewise affine space, and vector test functions to its
three-component counterpart.
\begin{subequations}\label{eq:ES-BDF3-MCF}
	\begin{align}
		\left(
		\frac{D_t^3\mathbf X^{m+1}}{\tau},
		\varphi^h\,\widetilde{\mathbf n}^{m+1}
		\right)^h_{\widetilde{\Gamma}^{m+1}}
		&+
		(1-\lambda^{m+1})
		\left(
		\mathcal H^{m+1},\varphi^h
		\right)^h_{\widetilde{\Gamma}^{m+1}}
		=0,
		\\
		\left(
		\mathcal H^{m+1},
		\widetilde{\mathbf n}^{m+1}\cdot\boldsymbol{\omega}^h
		\right)^h_{\widetilde{\Gamma}^{m+1}}
		&-
		\left(
		\nabla_{\widetilde{\Gamma}^{m+1}}\mathbf X^{m+1},
		\nabla_{\widetilde{\Gamma}^{m+1}}\boldsymbol{\omega}^h
		\right)_{\widetilde{\Gamma}^{m+1}}
		=0,
		\\
		W^{m+1}-W^m
        &+\tau\mathscr Q_m[d_{\mathrm M},\mathcal H]=0.
	\end{align}
\end{subequations}

The nonnegativity of $\mathscr Q_m[d_{\mathrm M},\mathcal H]$ again
gives $W^{m+1}\leq W^m$.

For surface diffusion, the normal velocity and surface-area dissipation are
\[
	V=\Delta_\Gamma\mathcal H,
	\qquad
	\frac{dW}{dt}
	=
	-\int_{\Gamma(t)}
	|\nabla_\Gamma\mathcal H|^2\,dA.
\]
The corresponding augmented velocity equation is
\[
	\partial_t\mathbf X\cdot\mathbf n
	=
	\Delta_\Gamma\mathcal H+\lambda(t)\mathcal H.
\]
Set $d_{\mathrm S,\Gamma}(u)=\|\nabla_\Gamma u\|_{L^2(\Gamma)}^2$.
The ES-BDF3 surface-diffusion system follows from
\eqref{eq:ES-BDF3-CDF} by replacing $\kappa$ with $\mathcal H$,
$\partial_s$ with $\nabla_\Gamma$, and $d_{\mathrm C}$ with $d_{\mathrm S}$,
using the surface finite element spaces above.

\section{Numerical experiments}\label{sec:numerics}

We examine geometric accuracy, energy decay, mesh distribution, and
nonlinear iteration counts for the numerical schemes presented in
Sections~\ref{sec:staggered} and~\ref{sec:extensions}.
Shrinking circles and spheres provide exact smooth benchmarks.

\subsection{Curve evolution}\label{Curve evolution experiments}

We consider CSF, curve diffusion (CDF), and anisotropic curve diffusion.
Since tangential motion may lead to different parametrizations of the same
geometric curve, we measure geometric differences using the manifold distance.
For two simple closed curves $\Gamma_1$ and $\Gamma_2$ enclosing bounded regions
$\Omega_1$ and $\Omega_2$, respectively, the manifold distance is defined by~\cite{jiang2024second}
\[
	\mathcal{M}\left(\Gamma_1,\Gamma_2\right)
	:=
	\left|
	\left(\Omega_1\backslash\Omega_2\right)
	\cup
	\left(\Omega_2\backslash\Omega_1\right)
	\right|
	=
	|\Omega_1|+|\Omega_2|-2|\Omega_1\cap\Omega_2|.
\]

When no exact solution is available, we use the Cauchy difference
\[
	\mathcal{E}_{\mathcal M,i}
	:=
	\mathcal{E}_{\mathcal M}(T,\tau_i,\tau_{i+1})
	=
	\mathcal{M}
	\left(
	\Gamma^{T/\tau_i},
	\Gamma^{T/\tau_{i+1}}
	\right).
\]
Under an asymptotic expansion in which the successive geometric
differences are proportional to $|\tau_i^p-\tau_{i+1}^p|$, an observed
order $p_i>0$ can be estimated, when a positive solution exists, from
\begin{equation}\label{eq:generalized-cauchy-order-manifold}
	\frac{\mathcal{E}_{\mathcal M,i}}
	{\mathcal{E}_{\mathcal M,i+1}}
	=
	\frac{\tau_i^{p_i}-\tau_{i+1}^{p_i}}
	{\tau_{i+1}^{p_i}-\tau_{i+2}^{p_i}}.
\end{equation}
For positive differences and $\tau_i>\tau_{i+1}>\tau_{i+2}>0$,
a positive root exists uniquely precisely when
\[
\frac{\mathcal E_{\mathcal M,i}}{\mathcal E_{\mathcal M,i+1}}>
\frac{\log(\tau_i/\tau_{i+1})}{\log(\tau_{i+1}/\tau_{i+2})}.
\]
The right-hand side of \eqref{eq:generalized-cauchy-order-manifold}
is strictly increasing in $p>0$, with the displayed threshold as its
$p\downarrow0$ limit. Zero differences or failure of this condition make
the positive-order estimate undefined.
For a constant refinement ratio $r=\tau_i/\tau_{i+1}>1$, this reduces to
\[
	p_i
	=
	\frac{
	\log(\mathcal E_{\mathcal M,i}/\mathcal E_{\mathcal M,i+1})
	}{
	\log(\tau_i/\tau_{i+1})
	}.
\]

In the experiments below, the spatial and temporal resolutions are generally
refined simultaneously. The resulting values of $p_i$ should therefore be
interpreted as observed orders along the specified coupled space--time
refinement paths, rather than as estimates obtained with a fixed spatial
discretization. In particular, along $\tau=Ch^{2/k}$, a second-order
spatial error is also $O(\tau^k)$; a slope close to $k$ on this path does
not isolate the temporal order. The symmetric-difference area measures
geometric agreement, not an error in a particular parametrization.

At the continuous level, the Lagrange multiplier vanishes for the flows
considered here under the nondegeneracy conditions stated in Section~\ref{sec:extensions}.
We therefore use
\[
	\mathcal{E}_{\lambda}(T,\tau_i)
	=
	\left|\lambda^{T/\tau_i}\right|,
	\qquad
	p_{\lambda,i}
	=
	\frac{
	\log\left(
	\mathcal{E}_{\lambda}(T,\tau_i)/
	\mathcal{E}_{\lambda}(T,\tau_{i+1})
	\right)}
	{\log(\tau_i/\tau_{i+1})}
\]
as a diagnostic of the discrete multiplier and its decay under refinement.

We also monitor the normalized energy $W(t)/W(0)$, the number of Newton
iterations, and the mesh ratio
\[
	\Psi^m
	=
	\frac{\max_{1\le j\le N}|\mathbf h_j^m|}
	{\min_{1\le j\le N}|\mathbf h_j^m|}.
\]
These quantities respectively describe the discrete energy evolution,
nonlinear iteration behavior, and mesh distribution.

\paragraph{Example 1 (Curve-shortening flow, CSF)}

The initial polygon has equally spaced vertices on the unit circle. The exact
solution is
\[
	\mathbf{X}_{\mathrm{ex}}(\rho,t)
	=
	\sqrt{1-2t}
	(\cos(2\pi\rho),\sin(2\pi\rho)),
	\quad
	\kappa_{\mathrm{ex}}
	=
	\frac{1}{\sqrt{1-2t}},
	\quad
	\rho\in\mathbb I,\quad t\in[0,0.5).
\]

For this example, define the geometric and curvature errors by
\[
\begin{aligned}
\mathcal E_{\mathcal M}(T,\tau)
 &=\mathcal M(\Gamma^{T/\tau},\Gamma_{\mathrm{ex}}(T)),\\
\mathcal E_\kappa(T,\tau)
 &=\|\kappa^{T/\tau}-\kappa_{\mathrm{ex}}(T)\|_\infty.
\end{aligned}
\]
For either error $\mathcal E$, the observed order between two resolutions is
\[
\operatorname{ord}(\tau_1,\tau_2)
=\frac{\log(\mathcal E(T,\tau_1)/\mathcal E(T,\tau_2))}
       {\log(\tau_1/\tau_2)}.
\]

\begin{figure}[htbp]
	\centering
	\includegraphics[width=0.95\linewidth]{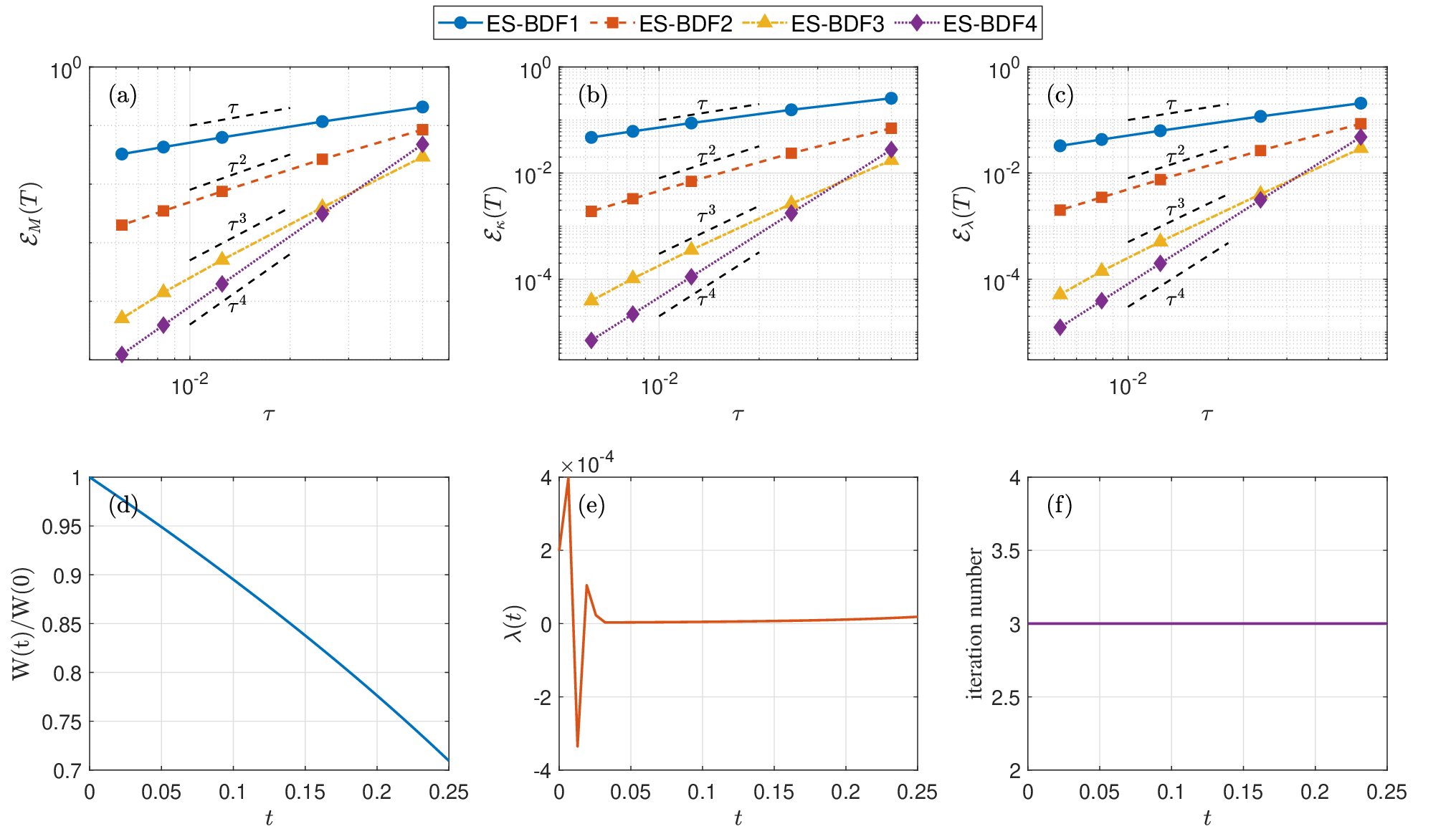}
	\caption{	Convergence results and representative diagnostics for CSF with a
	unit-circle initial curve. (a)--(c): log-log plots of
	$\mathcal{E}_{\mathcal M}$, $\mathcal{E}_{\kappa}$, and
	$\mathcal{E}_{\lambda}$ at $T=0.25$ for the ES-BDF\(k\) schemes;
	(d)--(f): evolution of the normalized perimeter, Lagrange multiplier,
	and Newton iteration count for ES-BDF3 with $N=1280$ and
	$\tau=1/160$.}
	\label{fig:CSF convergence order}
\end{figure}

Fig.~\ref{fig:CSF convergence order} shows the convergence results for the
shrinking-circle benchmark. Under the refinements used in this experiment,
the manifold-distance and curvature errors display the expected second-,
third-, and fourth-order trends for ES-BDF2, ES-BDF3, and ES-BDF4,
respectively. The multiplier magnitude
$\mathcal E_\lambda=|\lambda^{T/\tau}|$ also decreases at a high-order rate
under refinement. The second row of Fig.~\ref{fig:CSF convergence order} gives the numerical performance of the scheme ES-BDF3. 
As shown in this figure, the perimeter decreases monotonically, consistent with
Proposition~\ref{prop:energy}; the multiplier remains close to zero,
and the Newton iteration count always stays $3$ during the evolution.

\paragraph{Example 2 (Curve diffusion flow, CDF)}

We next consider CDF starting from a $4\times1$ ellipse. Since no exact
solution is available, convergence is assessed by Cauchy differences in the
manifold distance along coupled space--time refinement paths. We set
$T=0.25$ and use
$\tau=0.05h$, $\tau=0.05h^{2/3}$, and $\tau=0.05h^{1/2}$ for the
second-, third-, and fourth-order schemes, respectively.

\begin{figure}[htbp]
	\centering
	\includegraphics[width=0.95\linewidth]{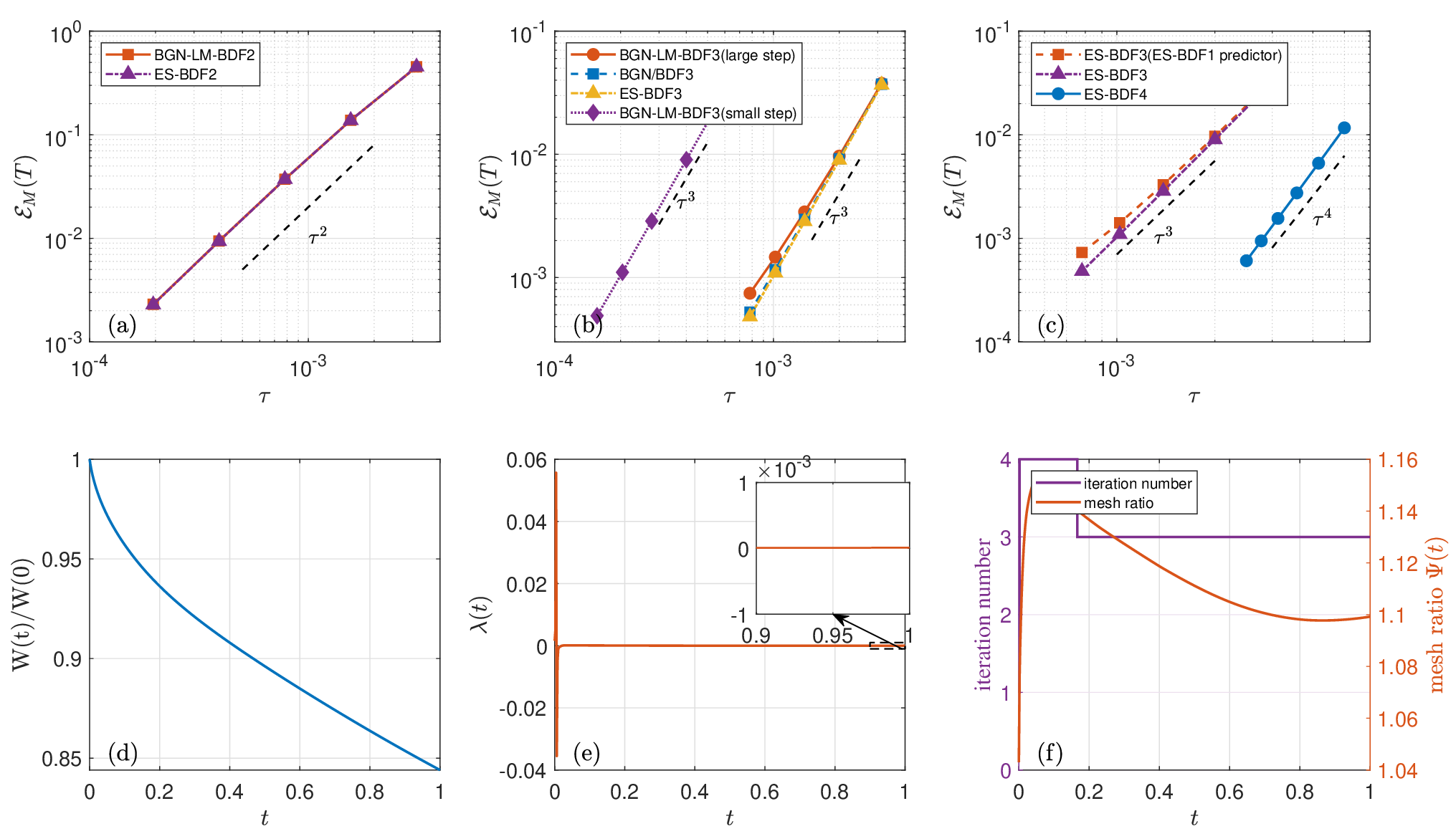}
	\caption{	Convergence results and diagnostics for CDF starting from a
	$4\times1$ ellipse. First row: Cauchy-type convergence tests for the
	reported BDF\(k\) schemes at $T=0.25$; Second row: ES-BDF3 diagnostics
	with $N=125$ and $\tau=1/500$ up to $T=1$: (d) normalized perimeter;
	(e) Lagrange multiplier; (f) Newton iteration count and mesh ratio
	$\Psi(t)$.}
	\label{fig:CDF}
\end{figure}

The first row of Fig.~\ref{fig:CDF} compares the observed convergence
behavior. ES-BDF2 and BGN-LM-BDF2 exhibit comparable second-order trends.
For the third-order schemes, ES-BDF3 and BGN/BDF3 display rates close to
three along the reported refinement path. BGN-LM-BDF3 shows a lower observed
rate on the coarser computations and approaches third-order behavior after
the time-step coefficient is reduced from $0.05$ to $0.01$.

The same experiment also compares two treatments of the midpoint data in
ES-BDF3. Replacing the staggered half-integer solve by a first-order
half-step predictor produces larger errors over the
tested refinement range. This numerical comparison supports the use of a
higher-order midpoint reconstruction in the Simpson energy update.
ES-BDF4 displays an approximately fourth-order trend along its reported
refinement path. The second row of Fig.~\ref{fig:CDF} shows representative ES-BDF3
diagnostics. The normalized perimeter decreases monotonically, the multiplier
remains close to zero, and the Newton iteration count and mesh ratio remain
bounded over the reported time interval.

\paragraph{Example 3 (Anisotropic curve diffusion)}

We next consider anisotropic curve diffusion using the formulation introduced
in Section~\ref{sec:extensions}. The initial curve is an ellipse with semi-major axis $2$ and
semi-minor axis $1$, and the anisotropic energy density is
\[
	\gamma(\theta)=1+\beta\cos(4\theta).
\]

Since no exact solution is available, convergence is assessed by the
generalized Cauchy procedure in
\eqref{eq:generalized-cauchy-order-manifold}. We take $\beta=0.05$ and use
the coupled refinement paths
$\tau=0.05h$, $\tau=0.05h^{2/3}$, and $\tau=0.05h^{1/2}$ for ES-BDF2,
ES-BDF3, and ES-BDF4, respectively. The convergence test is terminated at
$T=0.1$. Here $\gamma+\gamma''=1-15\beta\cos(4\theta)\geq0.25>0$,
so the test lies in the weakly anisotropic regime.

Table~\ref{tab:ASDF_cauchy} shows observed orders close to two and three for
ES-BDF2 and ES-BDF3 over most of the reported refinement range. The ES-BDF4
results exhibit observed orders between approximately $3.5$ and $3.7$ over
the tested range. The data therefore indicate convergence at a rate above three, but
do not yet resolve a fourth-order asymptotic regime.

\begin{table}[htbp]
	\centering
	\small
	\setlength{\tabcolsep}{4pt}
	\renewcommand{\arraystretch}{1.08}
	\begin{tabular}{ccc@{\qquad}ccc@{\qquad}ccc}
		\hline
		\multicolumn{3}{c}{ES-BDF2}
		&
		\multicolumn{3}{c}{ES-BDF3}
		&
		\multicolumn{3}{c}{ES-BDF4}
		\\
		$\tau$ & $\mathcal E_{\mathcal M}$ & order
		&
		$\tau$ & $\mathcal E_{\mathcal M}$ & order
		&
		$\tau$ & $\mathcal E_{\mathcal M}$ & order
		\\
		\hline
		$1/320$  & $1.555{\rm e}{-1}$ & --   & $1/180$ & $5.390{\rm e}{-2}$ & --   & $1/200$ & $3.071{\rm e}{-3}$ & -- \\
		$1/640$  & $3.621{\rm e}{-2}$ & 2.10 & $1/320$ & $8.684{\rm e}{-3}$ & 2.98 & $1/240$ & $1.382{\rm e}{-3}$ & 3.72 \\
		$1/1280$ & $8.575{\rm e}{-3}$ & 2.08 & $1/500$ & $2.083{\rm e}{-3}$ & 2.96 & $1/280$ & $6.994{\rm e}{-4}$ & 3.71 \\
		$1/2560$ & $1.947{\rm e}{-3}$ & 2.14 & $1/720$ & $6.282{\rm e}{-4}$ & 3.02 & $1/320$ & $3.920{\rm e}{-4}$ & 3.59 \\
		$1/5120$ & $3.930{\rm e}{-4}$ & 2.31 & $1/980$ & $2.465{\rm e}{-4}$ & 2.73 & $1/360$ & $2.364{\rm e}{-4}$ & 3.52 \\
		\hline
	\end{tabular}
	\caption{Convergence results for anisotropic curve diffusion with
	$\gamma(\theta)=1+\beta\cos(4\theta)$, $\beta=0.05$, and $T=0.1$. The coupled
	space--time refinement paths are $\tau=0.05h$ for ES-BDF2,
    $\tau=0.05h^{2/3}$ for ES-BDF3, and $\tau=0.05h^{1/2}$ for ES-BDF4.}
	\label{tab:ASDF_cauchy}
\end{table}

\begin{figure}[htbp]
	\centering
	\includegraphics[width=0.95\linewidth]{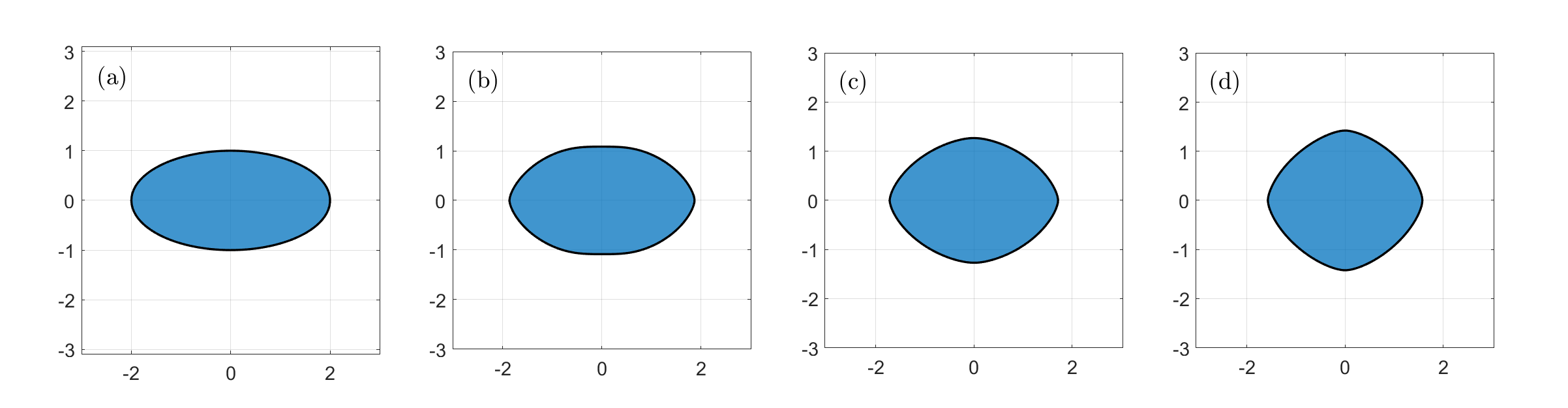}
	\caption{	Evolution of the curve under anisotropic curve diffusion at
	$t=0$, $0.1$, $0.25$, and $0.5$.
	}
	\label{aniSDF snapshots}
\end{figure}

\begin{figure}[htbp]
	\centering
	\includegraphics[width=0.72\linewidth]{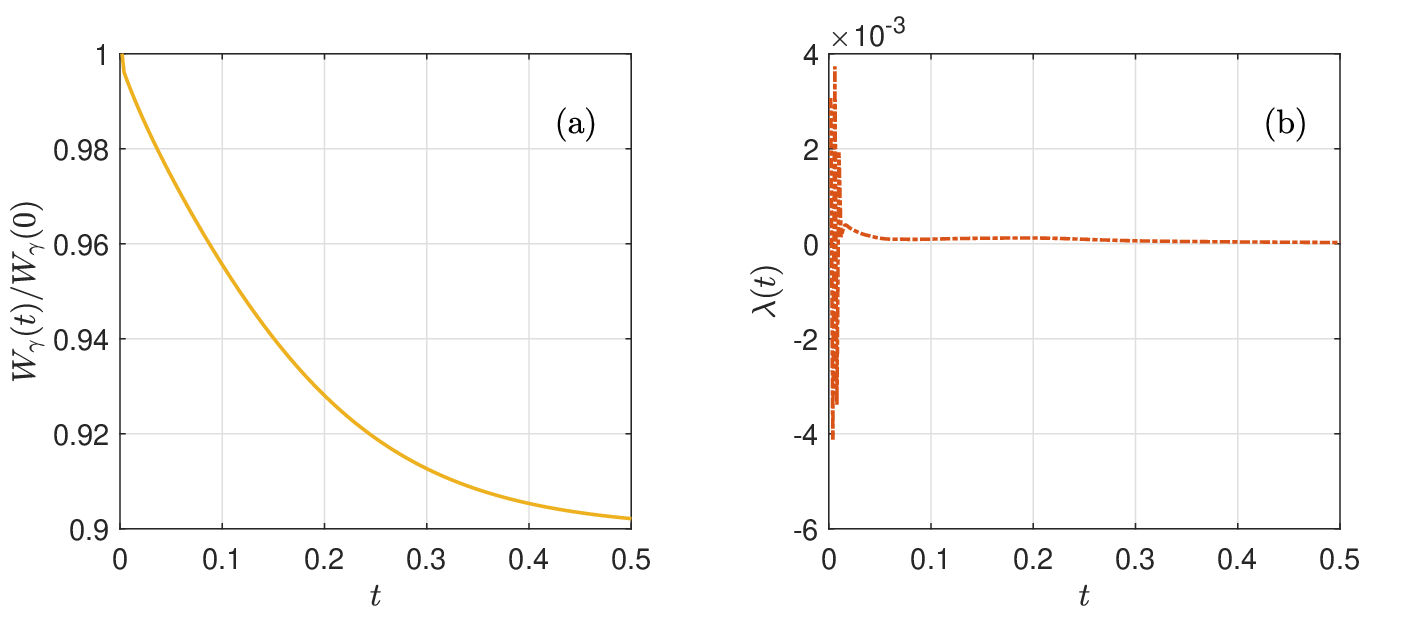}
	\caption{	Evolution of the normalized anisotropic energy (left) and Lagrange
	multiplier (right) up to $T=0.5$.
	}
	\label{fig:aniSDF_energy_lambda}
\end{figure}

For the longer simulation up to $T=0.5$, we use $N=125$ and
$\tau=0.05h^{2/3}$. Fig.~\ref{aniSDF snapshots} shows the evolution of the
curve under the anisotropic flow. Fig.~\ref{fig:aniSDF_energy_lambda}
reports the normalized anisotropic energy and the Lagrange multiplier.
The discrete anisotropic energy decreases monotonically, while the multiplier
remains close to zero over the reported interval, consistent with the
continuous relation.

\subsection{Surface evolution}\label{Surface evolution experiments}

To monitor the evolving triangulation, we use the two mesh-distribution
indicators introduced in \cite{jiang2024stable}. At $t=t_m$, they are
defined by
\begin{equation}\label{eq:surface-mesh-quality}
	r_h(t_m)
	:=
	\frac{
		\displaystyle
		\max_{1\le j\le J}
		\max\left\{
			\|\mathbf q_{j1}^m-\mathbf q_{j2}^m\|,
			\|\mathbf q_{j2}^m-\mathbf q_{j3}^m\|,
			\|\mathbf q_{j3}^m-\mathbf q_{j1}^m\|
		\right\}
	}{
		\displaystyle
		\min_{1\le j\le J}
		\min\left\{
			\|\mathbf q_{j1}^m-\mathbf q_{j2}^m\|,
			\|\mathbf q_{j2}^m-\mathbf q_{j3}^m\|,
			\|\mathbf q_{j3}^m-\mathbf q_{j1}^m\|
		\right\}
	},
\end{equation}
and
\begin{equation}
	r_a(t_m)
	:=
	\frac{
		\displaystyle\max_{1\le j\le J}|\sigma_j^m|
	}{
		\displaystyle\min_{1\le j\le J}|\sigma_j^m|
	}.
\end{equation}
Thus $r_h$ is the ratio of the largest to the smallest edge length, while
$r_a$ is the corresponding ratio of triangle areas. These quantities provide
simple diagnostics of edge-length and element-area variation during the
evolution. 

\paragraph{Example 4 (Mean curvature flow, MCF)}

We first consider MCF starting from the unit sphere. The exact solution
remains spherical with radius
\[
	R(t)=\sqrt{1-4t},
	\qquad t\in[0,1/4).
\]
We measure the geometric error by the maximal radial deviation
\[
	\mathcal E_{\mathrm{rad}}(T,\tau)
	:=
	\max_{\mathbf q\in\Gamma^{T/\tau}}
	\bigl|\|\mathbf q\|_2-R(T)\bigr|,
	\qquad
	\text{Order}
	=
	\frac{
		\log\left(
		\mathcal E_{\mathrm{rad}}(T,\tau_1)/
		\mathcal E_{\mathrm{rad}}(T,\tau_2)
		\right)}
	{\log(\tau_1/\tau_2)}.
\]

Following \cite{jiang2024stable}, the initial triangulation contains
$J=93608$ triangles and $K=46806$ vertices. The fine spatial resolution is
used to reduce the contribution of the spatial discretization error in this
temporal-refinement test. Observed temporal slopes are informative only
over the range where the spatial error and solver error are smaller than
the time-discretization error. 

\begin{figure}[htbp]
	\centering
	\includegraphics[width=0.95\linewidth]{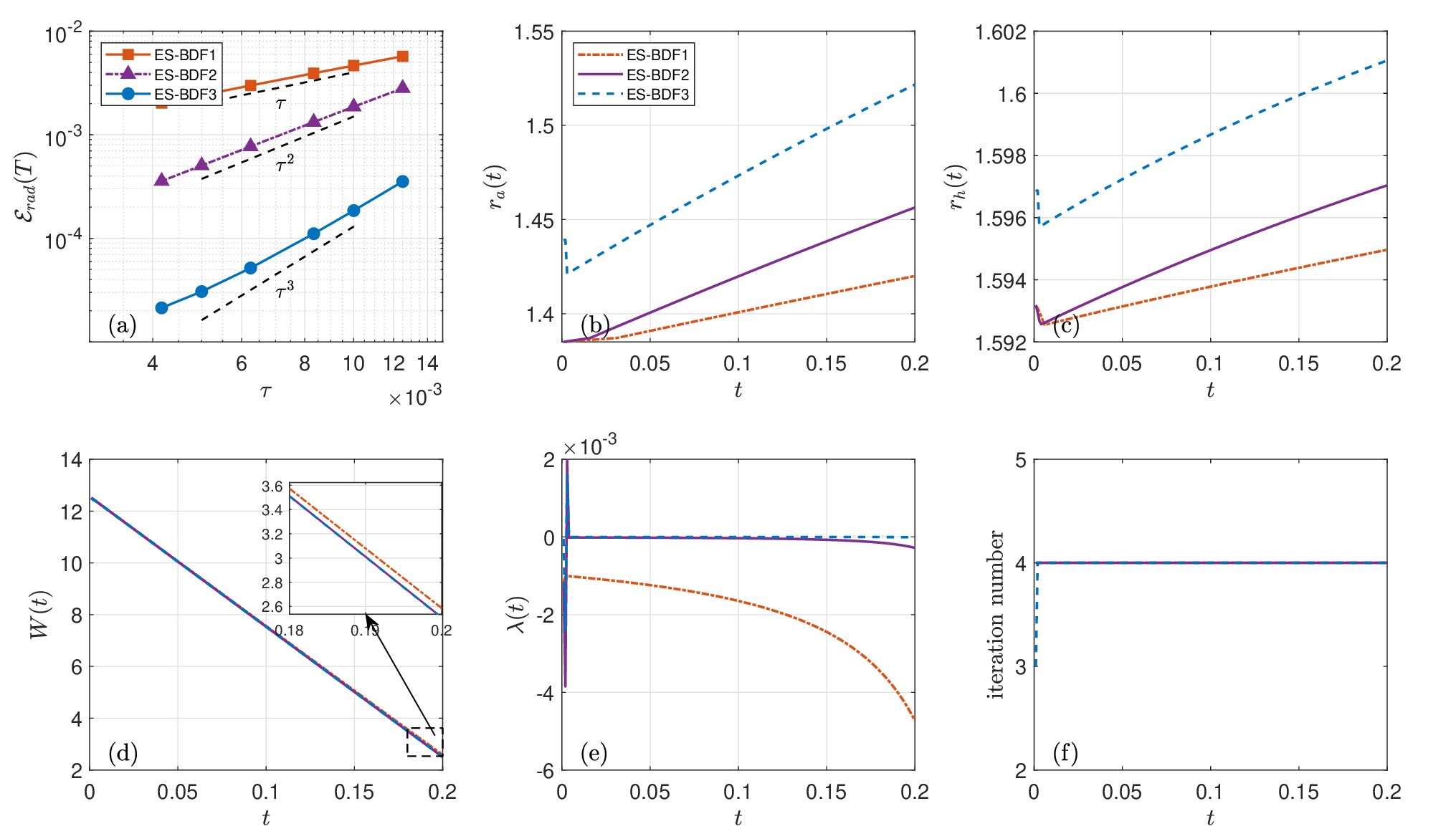}
	\caption{Results for MCF starting from the unit sphere, for ES-BDF1,
    ES-BDF2, and ES-BDF3. (a) Maximal radial deviation $\mathcal E_{\mathrm{rad}}$;
    (b) Triangle-area ratio $r_a(t)$. (c) Edge-length ratio $r_h(t)$;
    (d) Surface area $W(t)$;
    (e) Lagrange multiplier; (f) Newton iteration count.}
	\label{fig:MCF_sphere}
\end{figure}

For the three methods tested in Fig.~\ref{fig:MCF_sphere}(a),
ES-BDF1, ES-BDF2, and ES-BDF3 show approximately first-, second-, and
third-order trends, respectively. The surface area decreases
monotonically, and the mesh indicators
$r_h(t)$ and $r_a(t)$ remain bounded over the reported evolution.

Figs.~\ref{fig:MCF_sphere}(d)--(f) compare the surface-area, multiplier,
and Newton-iteration histories. The surface-area curves produced by the
three schemes are close. The Lagrange multipliers remain small, with smaller
magnitudes observed for the higher-order schemes in this test. The Newton
iteration counts are also of comparable size for the three discretizations.

\begin{figure}[htbp]
	\centering
	\includegraphics[width=0.95\linewidth]{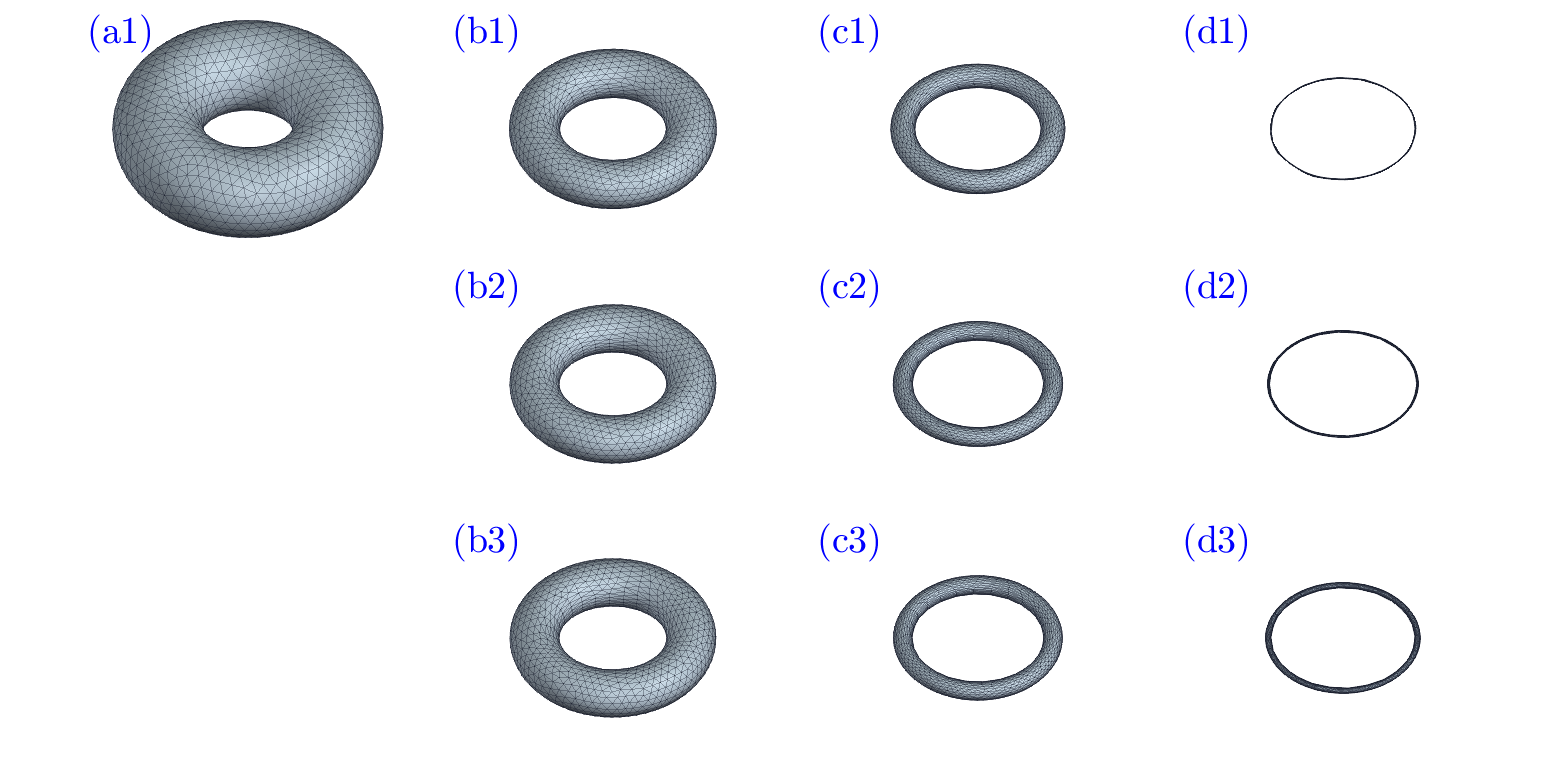}
	\caption{	Numerical evolution of the torus under MCF with
	$(J,K)=(2528,1264)$ and $\tau=10^{-3}$.
	First row: ES-BDF1 at $t=0$, $0.1$, $0.13$, and $0.145$;
	second row: ES-BDF2 at $t=0.1$, $0.13$, and $0.137$;
	third row: ES-BDF3 at $t=0.1$, $0.13$, and $0.135$.}
	\label{fig:MCF_torus_snapshots}
\end{figure}

We next consider the torus
\[
	(\sqrt{x_1^2+x_2^2}-R)^2+x_3^2=r^2,
	\qquad (R,r)=(1,0.5),
\]
using $(J,K)=(2528,1264)$ and $\tau=10^{-3}$.
Fig.~\ref{fig:MCF_torus_snapshots} shows the numerical evolution toward
the singular regime. The final pinch-off times for ES-BDF1, ES-BDF2, and
ES-BDF3 are $0.145$, $0.137$, and $0.135$, respectively, which are consistent with the
numerical value $0.138$ reported in \cite{barrett2008parametric1}.

\paragraph{Example 5 (Surface diffusion flow, SDF)}
Finally, we consider SDF of an $8\times1\times1$ cigar-like
surface obtained by capping a cylinder of diameter $1$ and length $7$ with
two hemispheres. The initial mesh contains $(J,K)=(2640,1322)$ elements and
vertices, respectively, and we choose $\tau=10^{-3}$. 

\begin{figure}[htbp]
    \includegraphics[width=0.95\linewidth]{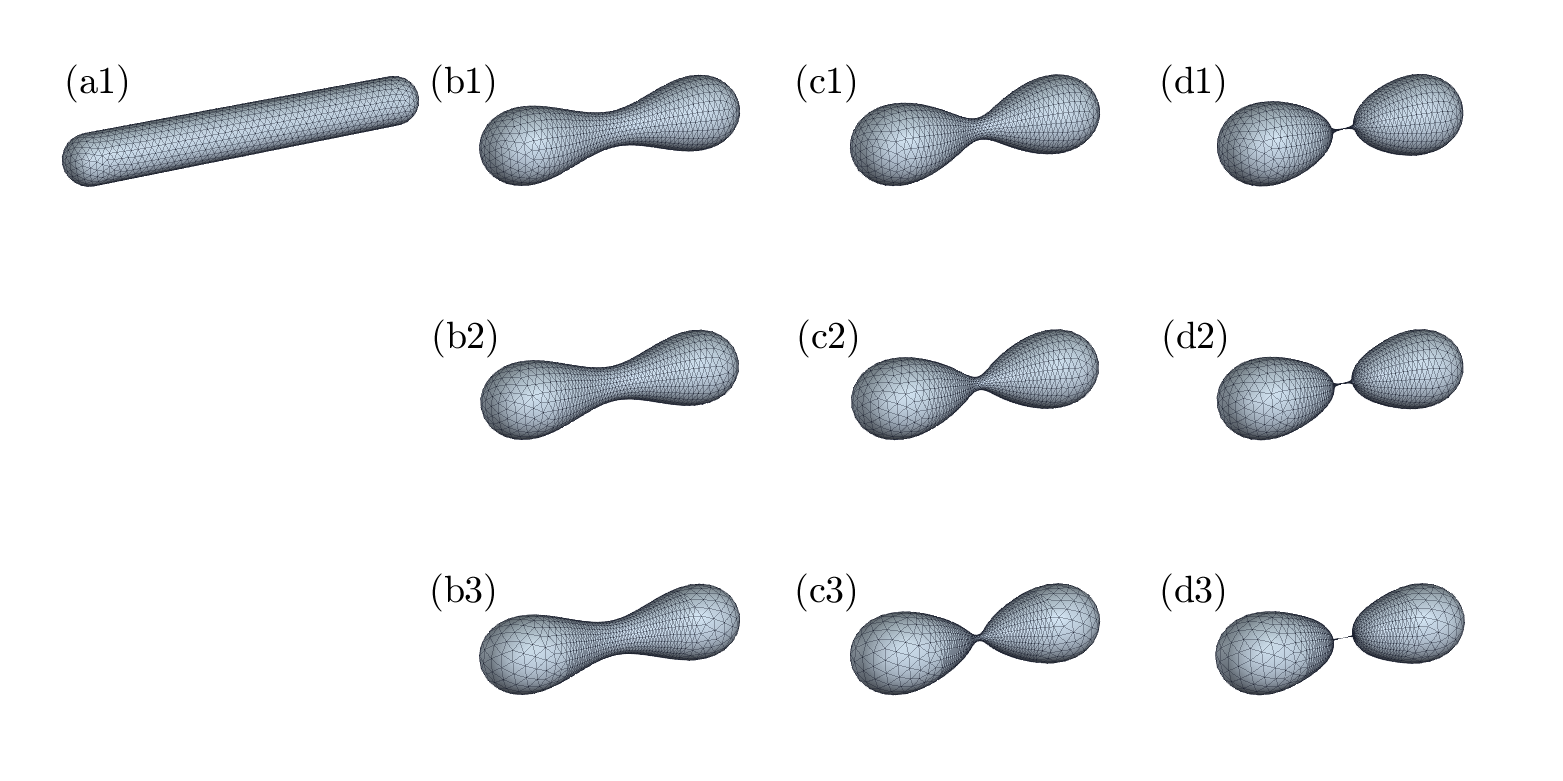}
	\caption{	Numerical evolution of the $8\times1\times1$ cigar-like surface under
	surface diffusion with $(J,K)=(2640,1322)$ and $\tau=10^{-3}$.
	First row: ES-BDF1 at $t=0$, $0.2$, $0.24$, and $0.257$;
	second row: ES-BDF2 at $t=0.2$, $0.24$, and $0.245$;
	third row: ES-BDF3 at $t=0.2$, $0.24$, and $0.243$.}
	\label{fig:SDF_cigar_snapshots}
\end{figure}

\begin{figure}[htbp]
	\includegraphics[width=0.95\linewidth]{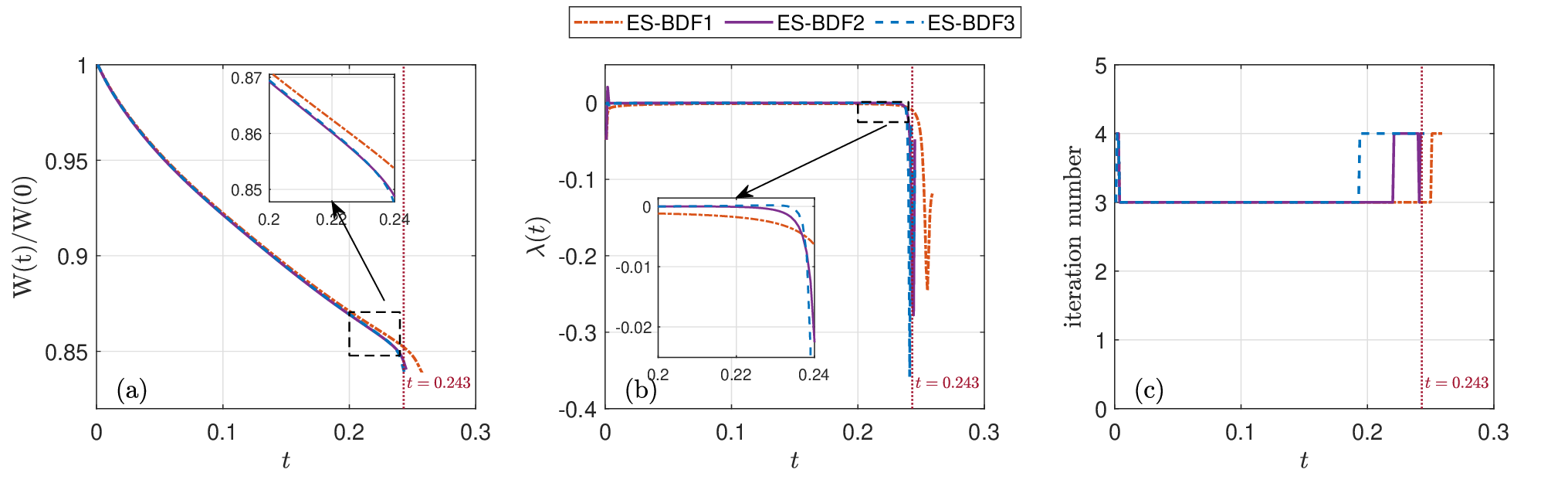}
	\caption{	Numerical diagnostics for ES-BDF\(k\), $k=1,2,3$, applied to surface
	diffusion of the cigar-like surface: (a) Normalized surface area;
	(b) Lagrange multiplier; (c) Newton iteration count.}
	\label{fig:SDF_cigar_diagnostics}
\end{figure}

Fig.~\ref{fig:SDF_cigar_snapshots} compares the evolutions produced by
ES-BDF1, ES-BDF2, and ES-BDF3. The surface develops a progressively narrower
neck and becomes nonconvex, and finally pinches off. The final pinch-off times are $0.257$, $0.245$,
and $0.243$ for ES-BDF1, ES-BDF2, and ES-BDF3, respectively. Compared to the pinch-off
resolution, the values produced by ES-BDF2 and ES-BDF3  are much closer than to
the value by ES-BDF1.

Fig.~\ref{fig:SDF_cigar_diagnostics} shows monotone area decay for
all three schemes. The ES-BDF2 and ES-BDF3 energy curves remain close,
while their separation from ES-BDF1 grows near the pinch-off time. The multiplier evolution
shows a similar pattern, and the Newton iteration counts are comparable.
Since the continuous Lagrange multiplier vanishes,
its discrete magnitude keeps very close to $0$ during the evolution, but exhibit
abruptly increasing when the evolution approaches the pinch-off time, shown in  
Fig.~\ref{fig:SDF_cigar_diagnostics}(b). The pinch-off time produced by ES-BDF3, $0.243$, is consistent with the numerical
value $0.2431$ reported in the literature~\cite{barrett2008parametric1}.

\section{Conclusions}\label{sec:conclusions}

We have developed energy-increment BDF methods for geometric flows within the BGN parametric finite element framework. The formulation replaces the multistep energy-derivative equation by a quadrature approximation of the integrated energy balance while retaining the BDF geometric equations. For BDF3 and BDF4, positive Simpson weights ensure decay of the original discrete geometric energy separately on each interleaved time grid for every nondegenerate solution satisfying the scalar constraint exactly. The analysis establishes conditional scalar consistency, a regularity criterion for the ES-BDF3 Newton system in curve shortening flow, and residual-based bounds on energy-balance errors. The smooth benchmarks exhibit accuracy trends consistent with the design orders under the reported refinements. The curve-diffusion comparison supports the importance of accurate midpoint data, while representative curve and surface tests document energy decay and mesh evolution.  A complete solvability and convergence analysis of the coupled scheme remains to be established, requiring joint control of initialization, predicted geometry, midpoint data, and nonlinear solution errors.

\bibliographystyle{siamplain}
	\bibliography{references}

@article{barrett2008parametric,
  doi = {10.1137/070700231},
  title={Parametric approximation of {Willmore} flow and related geometric evolution equations},
  author={Barrett, John W and Garcke, Harald and N{\"u}rnberg, Robert},
  journal={SIAM J. Sci. Comput.},
  volume={31},
  number={1},
  pages={225--253},
  year={2008},
  publisher={SIAM}
}

@article{barrett2008parametric1,
  doi = {10.1016/j.jcp.2007.11.023},
  title={On the parametric finite element approximation of evolving hypersurfaces in {${\mathbb R}^3$}},
  author={Barrett, John W and Garcke, Harald and N{\"u}rnberg, Robert},
  journal={J. Comput. Phys.},
  volume={227},
  number={9},
  pages={4281--4307},
  year={2008},
  publisher={Elsevier}
}

@article{jiang2024second,
  doi = {10.1016/j.jcp.2024.113220},
  title={A second-order in time, {BGN}-based parametric finite element method for geometric flows of curves},
  author={Jiang, Wei and Su, Chunmei and Zhang, Ganghui},
  journal={J. Comput. Phys.},
  volume={514},
  pages={113220},
  year={2024},
  publisher={Elsevier}
}

@article{barrett2007parametric,
  title={A parametric finite element method for fourth order geometric evolution equations},
  doi = {10.1016/j.jcp.2006.07.026},
  author={Barrett, John W and Garcke, Harald and N{\"u}rnberg, Robert},
  journal={J. Comput. Phys.},
  volume={222},
  number={1},
  pages={441--467},
  year={2007},
  publisher={Elsevier}
}

@incollection{barrett2020parametric,
  title={Parametric finite element approximations of curvature-driven interface evolutions},
  doi = {10.1016/bs.hna.2019.05.002},
  author={Barrett, John W and Garcke, Harald and N{\"u}rnberg, Robert},
  booktitle={Handbook of numerical analysis},
  volume={21},
  pages={275--423},
  year={2020},
  publisher={Elsevier}
}

@article{jiang2024stable,
  doi = {10.1137/23M1625597},
  title={Stable backward differentiation formula time discretization of {BGN}-based parametric finite element methods for geometric flows},
  author={Jiang, Wei and Su, Chunmei and Zhang, Ganghui},
  journal={SIAM J. Sci. Comput.},
  volume={46},
  number={5},
  pages={A2874--A2898},
  year={2024},
  publisher={SIAM}
}

@article{elliott2010modeling,
  title={Modeling and computation of two phase geometric biomembranes using surface finite elements},
  author={Elliott, Charles M and Stinner, Bj{\"o}rn},
  journal={J. Comput. Phys.},
  volume={229},
  number={18},
  pages={6585--6612},
  year={2010},
  publisher={Elsevier}
}

@article{garcke2025structure,
  doi = {10.1137/24M1687546},
  title={Structure-Preserving Parametric Finite Element Method for Surface Diffusion Based on {Lagrange} Multiplier Approaches},
  author={Garcke, Harald and Jiang, Wei and Su, Chunmei and Zhang, Ganghui},
  journal={SIAM J. Sci. Comput.},
  volume={47},
  number={3},
  pages={A1983--A2011},
  year={2025},
  publisher={SIAM}
}

@article{duan2024new,
  title={New artificial tangential motions for parametric finite element approximation of surface evolution},
  doi = {10.1137/23M1551857},
  author={Duan, Beiping and Li, Buyang},
  journal={SIAM J. Sci. Comput.},
  volume={46},
  number={1},
  pages={A587--A608},
  year={2024},
  publisher={SIAM}
}

@article{li2021energy,
  doi = {10.1016/j.jcp.2021.110658},
  title={An energy-stable parametric finite element method for anisotropic surface diffusion},
  author={Li, Yifei and Bao, Weizhu},
  journal={J. Comput. Phys.},
  volume={446},
  pages={110658},
  year={2021},
  publisher={Elsevier}
}

@article{jiang2025predictor,
  doi = {10.1016/j.jcp.2025.113901},
  title={Predictor-corrector, {BGN}-based parametric finite element methods for surface diffusion},
  author={Jiang, Wei and Su, Chunmei and Zhang, Ganghui and Zhang, Lian},
  journal={J. Comput. Phys.},
  volume={530},
  pages={113901},
  year={2025},
  publisher={Elsevier}
}

@article{Bao2023Symmetrized,
author = {Bao, Weizhu and Jiang, Wei and Li, Yifei},
title = {A Symmetrized Parametric Finite Element Method for Anisotropic Surface Diffusion of Closed Curves},
doi = {10.1137/22M1500575},
journal = {SIAM J. Numer. Anal.},
volume = {61},
number = {2},
pages = {617-641},
year = {2023},
}

@article{hu2022evolving,
  doi = {10.1007/s00211-022-01309-9},
  title={Evolving finite element methods with an artificial tangential velocity for mean curvature flow and {Willmore} flow},
  author={Hu, Jiashun and Li, Buyang},
  journal={Numer. Math.},
  volume={152},
  number={1},
  pages={127--181},
  year={2022},
  publisher={Springer}
}

@article{bao2021structure,
  title={A structure-preserving parametric finite element method for surface diffusion},
  doi = {10.1137/21M1406751},
  author={Bao, Weizhu and Zhao, Quan},
  journal={SIAM J. Numer. Anal.},
  volume={59},
  number={5},
  pages={2775--2799},
  year={2021},
  publisher={SIAM}
}

@article{shen2019new,
  title={A new class of efficient and robust energy stable schemes for gradient flows},
  doi = {10.1137/17M1150153},
  author={Shen, Jie and Xu, Jie and Yang, Jiang},
  journal={SIAM Rev.},
  volume={61},
  number={3},
  pages={474--506},
  year={2019},
  publisher={SIAM}
}

@article{shen2018scalar,
  doi = {10.1016/j.jcp.2017.10.021},
  title={The scalar auxiliary variable ({SAV}) approach for gradient flows},
  author={Shen, Jie and Xu, Jie and Yang, Jiang},
  journal={J. Comput. Phys.},
  volume={353},
  pages={407--416},
  year={2018},
  publisher={Elsevier}
}

@article{gao2026energy,
  title={An energy-stable minimal deformation rate scheme for mean curvature flow and surface diffusion},
  doi = {10.1137/25M1753838},
  author={Gao, Guangwei and Garcke, Harald and Li, Buyang and Tang, Rong},
  journal={SIAM J. Sci. Comput.},
  volume={48},
  number={1},
  pages={A103--A131},
  year={2026},
  publisher={SIAM}
}

@article{huang2022new,
  doi = {10.1016/j.cma.2022.114718},
  title={A new class of implicit--explicit {BDFk} {SAV} schemes for general dissipative systems and their error analysis},
  author={Huang, Fukeng and Shen, Jie},
  journal={Comput. Methods Appl. Mech. Engrg.},
  volume={392},
  pages={114718},
  year={2022},
  publisher={Elsevier}
}

@article{fu2024energy,
  doi = {10.1090/mcom/4021},
  title={Energy diminishing implicit-explicit {Runge}--{Kutta} methods for gradient flows},
  author={Fu, Zhaohui and Tang, Tao and Yang, Jiang},
  journal={Math. Comp.},
  volume={93},
  number={350},
  pages={2745--2767},
  year={2024}
}

@article{fu2025higher,
  doi = {10.1007/s11425-024-2337-3},
  title={Higher-order energy-decreasing exponential time differencing {Runge}--{Kutta} methods for gradient flows},
  author={Fu, Zhaohui and Shen, Jie and Yang, Jiang},
  journal={Sci. China Math.},
  volume={68},
  number={7},
  pages={1727--1746},
  year={2025},
  publisher={Springer}
}

@article{guo2025structure,
  title={Structure-preserving parametric finite element methods for anisotropic surface diffusion flow with minimal deformation formulation},
  doi = {10.1016/j.cpc.2025.109620},
  author={Guo, Yihang and Li, Meng},
  journal={Comput. Phys. Commun.},
  volume={313},
  pages={109620},
  year={2025},
  publisher={Elsevier}
}

@article{bao2025unified,
  title={A unified structure-preserving parametric finite element method for anisotropic surface diffusion},
  doi = {10.1090/mcom/4022},
  author={Bao, Weizhu and Li, Yifei},
  journal={Math. Comp.},
  volume={94},
  number={355},
  pages={2113--2149},
  year={2025}
}

@article{davi1990motion,
  title={On the motion of a phase interface by surface diffusion},
  doi = {10.1007/BF00945835},
  author={Davi, Fabrizio and Gurtin, Morton E},
  journal={Z. Angew. Math. Phys.},
  volume={41},
  number={6},
  pages={782--811},
  year={1990},
  publisher={Springer}
}

@book{cao2003geometric,
  title={Geometric curve evolution and image processing},
  author={Cao, Fr{\'e}d{\'e}ric},
  year={2003},
  publisher={Springer Science \& Business Media}
}

@misc{zhang2026arbitrary,
      title={Arbitrary-order structure-preserving discretizations for geometric curvature flows}, 
      author={Ganghui Zhang and Boris D. Andrews and Patrick E. Farrell},
      year={2026},
      howpublished={arXiv:2605.20371},
      doi={10.48550/arXiv.2605.20371},
}

@article{ma2026highorder,
  title   = {High-order energy-stable {BGN} parametric finite element methods for geometric flows},
  author  = {Ma, Shu and Rao, Qiqi},
  journal = {arXiv preprint arXiv:2608.29877},
  year    = {2026},
  doi     = {10.48550/arXiv.2608.29877}
}

@article{duan2021highorder,
  title={High-Order Fully Discrete Energy Diminishing Evolving Surface Finite Element Methods for a Class of Geometric Curvature Flows},
  author={Duan, Beiping and Li, Buyang and Zhang, Zhimin},
  journal={Ann. Appl. Math.},
  volume={37},
  number={4},
  pages={405--436},
  year={2021},
  doi={10.4208/aam.OA-2021-0007}
}
\end{document}